# ON ADAPTIVE ESTIMATION OF LINEAR FUNCTIONALS[1]

By T. Tony Cai and Mark G. Low

*University of Pennsylvania*

Adaptive estimation of linear functionals over a collection of parameter spaces is considered. A between-class modulus of continuity, a geometric quantity, is shown to be instrumental in characterizing the degree of adaptability over two parameter spaces in the same way that the usual modulus of continuity captures the minimax difficulty of estimation over a single parameter space. A general construction of optimally adaptive estimators based on an ordered modulus of continuity is given. The results are complemented by several illustrative examples.

**1. Introduction.** Adaptive estimation of linear functionals occupies an important position in the theory of nonparametric function estimation. As a step toward the goal of adaptive estimation, attention is first focused on the more concrete goal of developing a minimax theory over a fixed parameter space which can, for example, specify the smoothness of the function. This theory is now well developed, particularly in the white noise with drift model

$$(1) \qquad dY(t) = f(t)\,dt + \frac{1}{\sqrt{n}}\,dW(t), \qquad -\frac{1}{2} \le t \le \frac{1}{2},$$

where $W(t)$ is a standard Brownian motion. This model arises as an approximation to many other nonparametric models such as those of density estimation, nonparametric regression and spectral estimation. See, for example, [1, 2, 19, 21].

Based on white noise data, Ibragimov and Hasminskii [15] constructed linear estimators with the smallest maximum mean squared error over convex symmetric parameter spaces. Donoho and Liu [9] and Donoho [8] extended this theory to general convex parameter spaces in terms of a modulus of

Received June 2002; revised September 2004.
[1]Supported in part by NSF Grant DMS-03-06576.

continuity,

(2) $$\omega(\varepsilon, \mathcal{F}) = \sup\{|Tg - Tf| : \|g - f\|_2 \leq \varepsilon; f \in \mathcal{F}, g \in \mathcal{F}\}.$$

Affine estimators play a fundamental role in this theory. For a convex function class $\mathcal{F}$ and linear functional $T$, set the minimax affine risk $R_A^*(n, \mathcal{F}) = \inf_{\hat{T}_{\text{affine}}} \sup_{f \in \mathcal{F}} E(\hat{T} - Tf)^2$ and the minimax risk $R_N^*(n, \mathcal{F}) = \inf_{\hat{T}} \sup_{f \in \mathcal{F}} E(\hat{T} - Tf)^2$. Donoho and Liu [9] and Donoho [8] have shown that

(3) $$\frac{1}{8}\omega^2\left(\frac{1}{\sqrt{n}}, \mathcal{F}\right) \leq R_N^*(n, \mathcal{F}) \leq R_A^*(n, \mathcal{F}) \leq \omega^2\left(\frac{1}{\sqrt{n}}, \mathcal{F}\right)$$

and that the modulus can be used to construct the optimal affine procedure.

A natural way to extend the minimax theory to an adaptation theory is to construct estimators which are simultaneously near minimax over a collection of smoothness classes. In general, however, this goal cannot be realized. Lepski [17] was the first to give examples which demonstrated that rate optimal adaptation over a collection of Lipschitz classes is not possible when estimating a function at a point. Efromovich and Low [14] showed that this phenomenon is true in general over a collection of nested symmetric sets where the minimax rates are algebraic of different orders. See also [16].

On the other hand, the goal of fully rate adaptive estimation of linear functionals can sometimes be realized. When the minimax rate over each parameter space is slower than any algebraic rate, Cai and Low [5] have given examples of nested symmetric sets where fully adaptive estimators can be constructed. In addition, when the sets are not symmetric, there are also examples where rate adaptive estimators can be constructed. Such is the case for estimating monotone functions where an estimator can adapt over Lipschitz classes. See [20]. Other recent results can be found in [10, 11, 12, 13, 18].

Although the above-mentioned examples show that there are cases where fully rate adaptive estimators exist and other cases where fully rate adaptive estimators do not exist, to date there is no general theory that characterizes exactly when adaptation is possible. The present paper provides a general adaptation theory for estimating linear functionals. We develop a geometric understanding of the adaptive estimation problem analogous to that given by Donoho [8] for minimax theory. This theory describes exactly when fully rate adaptive estimators exist, and when they do not exist, the theory provides a general construction of estimators with minimum adaptation cost.

This paper and its companion papers Cai and Low [6, 7] develop a coherent approach to minimax theory, adaptive estimation and the construction of adaptive confidence intervals. The theory relies on two geometric quantities—a between-class modulus of continuity and an ordered modulus of continuity. For a pair of parameter spaces $\mathcal{F}_1$ and $\mathcal{F}_2$, the ordered modulus of continuity is defined by

(4) $\qquad \omega(\varepsilon, \mathcal{F}_1, \mathcal{F}_2) = \sup\{Tg - Tf : \|g - f\|_2 \leq \varepsilon; f \in \mathcal{F}_1, g \in \mathcal{F}_2\}.$

The ordered modulus of continuity is instrumental in the construction of the adaptive estimators given in Sections 2, 4 and 5. It is a quantity derived from the geometry of the graph of the linear functional $T$ between $\mathcal{F}_1$ and $\mathcal{F}_2$. It is also convenient to define a between-class modulus of continuity $\omega_+(\varepsilon, \mathcal{F}_1, \mathcal{F}_2)$ by

(5) $\qquad \omega_+(\varepsilon, \mathcal{F}_1, \mathcal{F}_2) = \sup\{|Tg - Tf| : \|g - f\|_2 \leq \varepsilon; f \in \mathcal{F}_1, g \in \mathcal{F}_2\}.$

Clearly, $\omega_+(\varepsilon, \mathcal{F}_1, \mathcal{F}_2) = \max\{\omega(\varepsilon, \mathcal{F}_1, \mathcal{F}_2), \omega(\varepsilon, \mathcal{F}_2, \mathcal{F}_1)\}$. When $\mathcal{F}_1 = \mathcal{F}_2 = \mathcal{F}$, $\omega(\varepsilon, \mathcal{F}, \mathcal{F}) = \omega_+(\varepsilon, \mathcal{F}, \mathcal{F})$ is the usual modulus of continuity over $\mathcal{F}$ and will be denoted by $\omega(\varepsilon, \mathcal{F})$ as in (2). We show that the between-class modulus can be used to characterize when adaptation is possible. This modulus captures the degree of adaptability over two parameter spaces in the same way that the usual modulus of continuity captures the minimax difficulty of estimation over a single parameter space.

We begin in Section 2 with a complete treatment of adaptation over an arbitrary pair of convex parameter spaces and any linear functional. In particular, we do not assume that the parameter spaces are nested or symmetric. A general construction for an optimally adaptive estimator is given. The adaptive estimator is based on appropriate tests between the parameter spaces which rely on a general understanding of the possible tradeoffs of bias and variance using the ordered moduli of continuity.

The theory shows that there are three main cases in terms of the cost of adaptation. We shall call the first case the regular one where, as in the case of estimating a function at a point over Lipschitz classes, the cost of adaptation is a logarithmic factor of the noise level. In the second case, full adaptation is possible as in the examples considered in [5, 18]. More dramatically, in the third case, the cost of adaptation is much greater than in the regular case. The cost of adaptation in this case is a power of the noise level. Examples of all three cases are given in Section 3.

Section 2 gives a geometric characterization of adaptation and shows the fundamental role played by the between-class and ordered moduli of continuity in this theory. The adaptation theory over two spaces in turn provides a fundamental building block for adaptation over richer collections of parameter spaces. In Section 4 we extend this theory to any collection of finitely many nested convex spaces, and under mild regularity conditions on the modulus to finitely many nonnested convex parameter spaces. The focus of this section is on the construction of an estimator with minimum adaptation cost. In Section 5 we further generalize the results to a continuum of parameter spaces.

**2. Adaptation over two parameter spaces.** In this section we give a complete development of adaptation over an arbitrary pair of convex parameter spaces $\mathcal{F}_1$ and $\mathcal{F}_2$ with $\mathcal{F}_1 \cap \mathcal{F}_2 \neq \varnothing$ and any linear functional $T$. We first derive a benchmark for the performance over $\mathcal{F}_2$ of any minimax rate optimal estimator over $\mathcal{F}_1$. The benchmark is given in terms of the between-class modulus of continuity. A general construction for an optimally adaptive estimator is then given. The adaptive procedure is built on a test between the parameter spaces which is based on the tradeoffs of bias and variance using the ordered moduli of continuity. Taken together these results show that the moduli of continuity captures the degree to which adaptation is possible.

Throughout the paper, we denote by $C$ a generic constant that may vary from place to place.

2.1. *Lower bound on the cost of adaptation.* Let the ordered modulus of continuity $\omega(\varepsilon, \mathcal{F}_1, \mathcal{F}_2)$ be defined as in (4) and the between-class modulus be given as in (5). Note that $\omega(\varepsilon, \mathcal{F}_1, \mathcal{F}_2)$ does not necessarily equal $\omega(\varepsilon, \mathcal{F}_2, \mathcal{F}_1)$. It is however clear that the modulus $\omega(\varepsilon, \mathcal{F}_1, \mathcal{F}_2)$ is an increasing function of $\varepsilon$. Moreover, if $\mathcal{F}_1$ and $\mathcal{F}_2$ are convex with $\mathcal{F}_1 \cap \mathcal{F}_2 \neq \varnothing$, then for a linear functional $T$ the modulus $\omega(\varepsilon, \mathcal{F}_1, \mathcal{F}_2)$ is also a concave function of $\varepsilon$. See [6]. Note also that although $\omega_+$ need not be concave, it follows from the concavity of the ordered modulus of continuity that for $D \geq 1$,

(6) $$\omega_+(D\varepsilon, \mathcal{F}_1, \mathcal{F}_2) \leq D\omega_+(\varepsilon, \mathcal{F}_1, \mathcal{F}_2).$$

The following result gives the lower bound for the maximum risk over $\mathcal{F}_2$ for minimax rate optimal estimators over $\mathcal{F}_1$.

THEOREM 1. *Let $T$ be a linear functional and let $\mathcal{F}_1$ and $\mathcal{F}_2$ be parameter spaces with $\mathcal{F}_1 \cap \mathcal{F}_2 \neq \varnothing$ and $\omega(\varepsilon, \mathcal{F}_1) \leq \omega(\varepsilon, \mathcal{F}_2)$ for all sufficiently small $0 < \varepsilon \leq \varepsilon_0$. Suppose that $\hat{T}$ is an estimator of $Tf$ based on the white noise data* (1) *satisfying*

(7) $$\sup_{f \in \mathcal{F}_1} E_f(\hat{T} - Tf)^2 \leq c_*^2 \omega^2\left(\frac{1}{\sqrt{n}}, \mathcal{F}_1\right)$$

*for some constant $c_* > 0$. Let $\gamma_n = \max\{e, \frac{\omega_+(1/\sqrt{n}, \mathcal{F}_1, \mathcal{F}_2)}{c_* \omega(1/\sqrt{n}, \mathcal{F}_1)}\}$. Then there exists some fixed constant $c > 0$ such that, for all sufficiently large $n$,*

(8) $$\sup_{f \in \mathcal{F}_2} E_f(\hat{T} - Tf)^2 \geq c\left\{\omega_+^2\left(\sqrt{\frac{\ln \gamma_n}{n}}, \mathcal{F}_1, \mathcal{F}_2\right) + \omega^2\left(\frac{1}{\sqrt{n}}, \mathcal{F}_2\right)\right\}.$$

PROOF. We shall only consider the case where $\mathcal{F}_1$ and $\mathcal{F}_2$ are closed and norm bounded. The general case is proved by taking limits of this case as in Section 14 of [8].

For the case of $\gamma_n = e$, (8) follows directly from (3). Now assume that $\gamma_n > e$. Then $\sup_{f \in \mathcal{F}_1} E_f(\hat{T} - Tf)^2 \leq \gamma_n^{-2} \omega_+^2(\frac{1}{\sqrt{n}}, \mathcal{F}_1, \mathcal{F}_2)$. Choose $f_{1,n} \in \mathcal{F}_1$ and $f_{2,n} \in \mathcal{F}_2$ such that $\|f_{1,n} - f_{2,n}\|_2 \leq \sqrt{\frac{\ln \gamma_n}{n}}$ and such that the between-class modulus is attained at $\{f_{1,n}, f_{2,n}\} : |Tf_{2,n} - Tf_{1,n}| = \omega_+(\sqrt{\frac{\ln \gamma_n}{n}}, \mathcal{F}_1, \mathcal{F}_2)$. It then follows from the constrained risk inequality of Brown and Low [3] and equation (6) that

$$E_{f_{2,n}}(\hat{T} - Tf_{2,n})^2 \geq \left(\omega_+\left(\sqrt{\frac{\ln \gamma_n}{n}}, \mathcal{F}_1, \mathcal{F}_2\right) - \gamma_n^{-1/2}\omega_+\left(\frac{1}{\sqrt{n}}, \mathcal{F}_1, \mathcal{F}_2\right)\right)^2$$

$$\geq (1 - e^{-1/2})^2 \omega_+^2\left(\sqrt{\frac{\ln \gamma_n}{n}}, \mathcal{F}_1, \mathcal{F}_2\right),$$

and hence by equation (3) $\sup_{f \in \mathcal{F}_2} E_f(\hat{T} - Tf)^2 \geq \frac{1}{16}\{\omega_+^2(\sqrt{\frac{\ln \gamma_n}{n}}, \mathcal{F}_1, \mathcal{F}_2) + \omega^2(\frac{1}{\sqrt{n}}, \mathcal{F}_2)\}$. □

Theorem 1 considers the performance over $\mathcal{F}_2$ of estimators which are minimax rate optimal over $\mathcal{F}_1$. This is a particularly important case but we shall also need a more general bound when we discuss adaptation over collections of parameter spaces. The proof of the following theorem is similar to that of Theorem 1 and is thus omitted.

THEOREM 2. *Consider two function classes $\mathcal{F}_1$ and $\mathcal{F}_2$ with $\mathcal{F}_1 \cap \mathcal{F}_2 \neq \varnothing$. Let $T$ be a linear functional and suppose that*

$$(9) \qquad \sup_{f \in \mathcal{F}_1} E_f(\hat{T} - Tf)^2 \leq \gamma_n^{-2}\omega_+^2\left(\frac{1}{\sqrt{n}}, \mathcal{F}_1, \mathcal{F}_2\right)$$

*for some $\gamma_n > 1$. Then for any $0 < \rho \leq 1$,*

$$(10) \quad \sup_{f \in \mathcal{F}_2}[E_f(\hat{T} - Tf)^2]^{1/2}$$

$$\geq \omega_+\left(\sqrt{\frac{\rho \ln \gamma_n^2}{n}}, \mathcal{F}_1, \mathcal{F}_2\right) - \gamma_n^{-(1-\rho)}\omega_+\left(\frac{1}{\sqrt{n}}, \mathcal{F}_1, \mathcal{F}_2\right).$$

2.2. *Construction of optimally adaptive procedure.* We now turn to a general construction of an adaptive procedure for any given linear functional $T$ over any two convex parameter spaces $\mathcal{F}_1$ and $\mathcal{F}_2$ with nonempty intersection $\mathcal{F}_1 \cap \mathcal{F}_2 \neq \varnothing$.

Before describing the adaptive procedure first focus attention on each parameter space separately. If it were known that $f \in \mathcal{F}_i$ then the theory of Donoho and Liu yields linear estimators $\hat{T}_i$ which satisfy $\sup_{f \in \mathcal{F}_i} E(\hat{T}_i -$ $Tf)^2 \leq \omega^2(\frac{1}{\sqrt{n}}, \mathcal{F}_i)$. Moreover these estimators are minimax rate optimal over $\mathcal{F}_i$. The adaptive procedure is then based on a test between $\mathcal{F}_1$ and $\mathcal{F}_2$. If the test accepts $\mathcal{F}_1$ then the procedure uses $\hat{T}_1$ whereas if it rejects $\mathcal{F}_1$ the procedure uses a minimax rate optimal procedure over $\mathcal{F}_1 \cup \mathcal{F}_2$. The test is designed in such a way that if $f \in \mathcal{F}_1$ it has a small probability of rejecting $\mathcal{F}_1$. On the other hand if $f \in \mathcal{F}_2$ and the bias of $\hat{T}_1$ is large the test has only a small probability of accepting $\mathcal{F}_1$.

In the case where $\mathcal{F}_1$ and $\mathcal{F}_2$ are nonnested convex parameter spaces it is clear that an implementation of this approach requires a minimax analysis for sets which are not convex. The reason is that we need to know the minimax risk and minimax rate optimal procedure over the union $\mathcal{G} = \mathcal{F}_1 \cup \mathcal{F}_2$, which is in general nonconvex. Such a theory has been given [6] where it was shown that if $\mathcal{G}$ is a union of a finite number of closed convex parameter spaces, the minimax risk is of the order $\omega^2(\frac{1}{\sqrt{n}}, \mathcal{G})$. Moreover, explicit rate optimal procedures, say $\hat{T}_2^*$, were constructed which for $\mathcal{G} = \mathcal{F}_1 \cup \mathcal{F}_2$ satisfy

$$\sup_{f \in \mathcal{G}} E(\hat{T}_2^* - Tf)^4 \leq C\omega^4\left(\frac{1}{\sqrt{n}}, \mathcal{G}\right). \tag{11}$$

In the adaptive procedure $\hat{T}_2^*$ is used whenever $\mathcal{F}_1$ is rejected.

The test between $\mathcal{F}_1$ and $\mathcal{F}_2$ is based on a comparison of linear estimators which trade bias and variance over $\mathcal{F}_1$ and $\mathcal{F}_2$ in a precise way. This trading of bias and variance is based on results in [6] which show how to use the ordered modulus of continuity to construct a linear procedure which has upper bounds for the bias over one parameter space and lower bounds for the bias over the other parameter space. More specifically, for two convex sets $\mathcal{F}$ and $\mathcal{H}$ with $\mathcal{F} \cap \mathcal{H} \neq \varnothing$, a linear estimator $\hat{T}$ is given which has variance and bias satisfying

$$\text{Var}(\hat{T}) = E(\hat{T} - E\hat{T})^2 \leq V, \tag{12}$$

$$\sup_{f \in \mathcal{F}}(E\hat{T} - Tf) \leq \tfrac{1}{2}\sup_{\varepsilon > 0}(\omega(\varepsilon, \mathcal{F}, \mathcal{H}) - \sqrt{nV}\varepsilon), \tag{13}$$

$$\inf_{f \in \mathcal{H}}(E\hat{T} - Tf) \geq -\tfrac{1}{2}\sup_{\varepsilon > 0}(\omega(\varepsilon, \mathcal{F}, \mathcal{H}) - \sqrt{nV}\varepsilon). \tag{14}$$

For a given bound $V$ on the variance this theory leads to two linear estimators by interchanging the roles of $\mathcal{F}$ and $\mathcal{H}$. In our context we make two different choices for $V$. For $1 \leq i \neq j \leq 2$ let

$$\gamma_{i,j} = \max\left(e, \frac{\omega(1/\sqrt{n}, \mathcal{F}_i, \mathcal{F}_j)}{\omega(1/\sqrt{n}, \mathcal{F}_1)}\right) \quad \text{and}$$

$$\gamma_+ = \max(\gamma_{1,2}, \gamma_{2,1}) = \max\left(e, \frac{\omega_+(1/\sqrt{n}, \mathcal{F}_1, \mathcal{F}_2)}{\omega(1/\sqrt{n}, \mathcal{F}_1)}\right) \tag{15}$$

and set
$$\sigma_{i,j}^2 = \frac{1}{\ln \gamma_{i,j}} \omega^2\left(\sqrt{\frac{\ln \gamma_{i,j}}{n}}, \mathcal{F}_i, \mathcal{F}_j\right).$$

The estimators $\hat{T}_{i,j}$ for $1 \leq i \neq j \leq 2$ needed in the test are linear estimators satisfying (12)–(14) for $\mathcal{F} = \mathcal{F}_i$, $\mathcal{H} = \mathcal{F}_j$ and $V = \sigma_{i,j}^2$.

The test given below relies on an understanding of the bias properties of $\hat{T}_{1,2}$ and $\hat{T}_{2,1}$. Simple bounds on the bias are easy to obtain from (13) and (14) since

$$\sup_{\varepsilon > 0}(\omega(\varepsilon, \mathcal{F}_i, \mathcal{F}_j) - \varepsilon\sqrt{n\sigma_{i,j}^2})$$

$$= \sup_{\varepsilon \leq \sqrt{\ln \gamma_{i,j}/n}}\left(\omega(\varepsilon, \mathcal{F}_i, \mathcal{F}_j) - \varepsilon\sqrt{\frac{n}{\ln \gamma_{i,j}}} \omega\left(\sqrt{\frac{\ln \gamma_{i,j}}{n}}, \mathcal{F}_i, \mathcal{F}_j\right)\right)$$

$$\leq \omega\left(\sqrt{\frac{\ln \gamma_{i,j}}{n}}, \mathcal{F}_i, \mathcal{F}_j\right).$$

The test is based on a comparison of the estimator $\hat{T}_1$ and both $\hat{T}_{1,2}$ and $\hat{T}_{2,1}$.

Note that if $f \in \mathcal{F}_1$,

(16)
$$E(\hat{T}_1 - \hat{T}_{1,2}) = E(\hat{T}_1 - Tf) - E(\hat{T}_{1,2} - Tf)$$
$$\geq -\omega\left(\frac{1}{\sqrt{n}}, \mathcal{F}_1\right) - \omega\left(\sqrt{\frac{\ln \gamma_{1,2}}{n}}, \mathcal{F}_1, \mathcal{F}_2\right) = -b_{1,2},$$

(17)
$$E(\hat{T}_1 - \hat{T}_{2,1}) = E(\hat{T}_1 - Tf) - E(\hat{T}_{2,1} - Tf)$$
$$\leq \omega\left(\frac{1}{\sqrt{n}}, \mathcal{F}_1\right) + \omega\left(\sqrt{\frac{\ln \gamma_{2,1}}{n}}, \mathcal{F}_2, \mathcal{F}_1\right) = b_{2,1},$$

(18)
$$\operatorname{Var}(\hat{T}_1 - \hat{T}_{1,2}) \leq 2\left(\omega^2\left(\frac{1}{\sqrt{n}}, \mathcal{F}_1\right) + \frac{1}{\ln \gamma_{1,2}} \omega^2\left(\sqrt{\frac{\ln \gamma_{1,2}}{n}}, \mathcal{F}_1, \mathcal{F}_2\right)\right)$$
$$= v_{1,2},$$

(19)
$$\operatorname{Var}(\hat{T}_1 - \hat{T}_{2,1}) \leq 2\left(\omega^2\left(\frac{1}{\sqrt{n}}, \mathcal{F}_1\right) + \frac{1}{\ln \gamma_{2,1}} \omega^2\left(\sqrt{\frac{\ln \gamma_{2,1}}{n}}, \mathcal{F}_2, \mathcal{F}_1\right)\right)$$
$$= v_{2,1}.$$

Hence, if $f \in \mathcal{F}_1$ it is easy to select a value so that the chance that $\hat{T}_1 - \hat{T}_{2,1}$ is greater than that value is small. Likewise it is easy to select another value so that the chance that $\hat{T}_1 - \hat{T}_{1,2}$ is less than that value is small. A careful selection of these values leads to the following test between $\mathcal{F}_1$ and $\mathcal{F}_2$:

$$(20) \quad I_n = 1\left(\hat{T}_{1,2} - 5b_{1,2} - 4\omega\left(\frac{1}{\sqrt{n}}, \mathcal{G}\right) \leq \hat{T}_1 \leq \hat{T}_{2,1} + 5b_{2,1} + 4\omega\left(\frac{1}{\sqrt{n}}, \mathcal{G}\right)\right).$$

The value $I_n = 1$ corresponds to accepting $\mathcal{F}_1$, in which case $\hat{T}_1$ is used. The value $I_n = 0$ corresponds to rejecting $\mathcal{F}_1$, in which case $\hat{T}_2^*$ is used. The adaptive estimator can then be written as

$$(21) \quad \hat{T} = I_n \hat{T}_1 + (1 - I_n)\hat{T}_2^*,$$

where $\hat{T}_1$ satisfies $\sup_{\mathcal{F}_1} E(\hat{T}_1 - Tf)^2 \leq \omega^2(\frac{1}{\sqrt{n}}, \mathcal{F}_1)$ and $\hat{T}_2^*$ satisfies (11).

2.3. *Adaptivity of the procedure.* In the previous subsection an estimator $\hat{T}$ was constructed based on a test between two parameter spaces. In this section we show that this estimator is adaptively rate optimal over $\mathcal{F}_1$ and $\mathcal{F}_2$. As a consequence it is shown that the lower bound for adaptation between $\mathcal{F}_1$ and $\mathcal{F}_2$ as given in Theorem 1 is sharp. The following theorem summarizes these results.

THEOREM 3. *Suppose $\mathcal{F}_1$ and $\mathcal{F}_2$ are two closed convex parameter spaces with $\mathcal{F}_1 \cap \mathcal{F}_2 \neq \varnothing$ and $\omega(\varepsilon, \mathcal{F}_1) \leq \omega(\varepsilon, \mathcal{F}_2)$. The estimator $\hat{T}$ defined in (21) satisfies for some fixed $C > 0$*

$$(22) \quad \sup_{f \in \mathcal{F}_1} E(\hat{T} - Tf)^2 \leq C\omega^2\left(\frac{1}{\sqrt{n}}, \mathcal{F}_1\right)$$

*and*

$$(23) \quad \sup_{f \in \mathcal{F}_2} E(\hat{T} - Tf)^2 \leq C\left\{\omega_+^2\left(\sqrt{\frac{\ln \gamma_+}{n}}, \mathcal{F}_1, \mathcal{F}_2\right) + \omega^2\left(\frac{1}{\sqrt{n}}, \mathcal{F}_2\right)\right\},$$

*where $\gamma_+$ is defined in (15).*

In light of the lower and upper bounds given in Theorems 1 and 3 we give the following definition.

DEFINITION 1. We shall call an estimator $\hat{T}$ *optimally adaptive* over $\mathcal{F}_1$ and $\mathcal{F}_2$ if it satisfies both (22) and (23).

REMARK. The estimator $\hat{T}$ defined in (21) is also adaptive between $\mathcal{F}_1$ and $\mathcal{G} = \mathcal{F}_1 \cup \mathcal{F}_2$. Note that (23) is equivalent to

$$(24) \quad \sup_{f \in \mathcal{G}} E(\hat{T} - Tf)^2 \leq C\left(\omega_+^2\left(\sqrt{\frac{\ln \gamma_+^*}{n}}, \mathcal{F}_1, \mathcal{G}\right) + \omega^2\left(\frac{1}{\sqrt{n}}, \mathcal{G}\right)\right),$$

where $\gamma_+^* = \max(e, \frac{\omega_+(1/\sqrt{n}, \mathcal{F}_1, \mathcal{G})}{\omega(1/\sqrt{n}, \mathcal{F}_1)})$. Therefore $\hat{T}$ attains the exact minimax rate of convergence over $\mathcal{F}_1$ and attains the lower bound on adaptation over $\mathcal{G}$ as given in Theorem 1.

As mentioned in the previous section, the estimator $\hat{T}$ was constructed by testing between $\mathcal{F}_1$ and $\mathcal{F}_2$. The proof of Theorem 3 is based on a precise analysis of the properties of the test as given in the following lemmas.

This test is constructed so that for $f \in \mathcal{F}_1$ the probability of rejecting $\mathcal{F}_1$ is small. Lemma 1 below provides a specific bound on the rejection of $\mathcal{F}_1$ when $f \in \mathcal{F}_1$.

LEMMA 1. *If $f \in \mathcal{F}_1$, then*

$$P(I_n = 0) \leq \frac{\omega^4(1/\sqrt{n}, \mathcal{F}_1)}{\omega^4(1/\sqrt{n}, \mathcal{G})}. \tag{25}$$

PROOF. First note that for a standard normal random variable $Z$, $P(Z \geq \lambda) \leq \exp(-\frac{\lambda^2}{2})$ holds for all $\lambda \geq 0$. It then follows from (16)–(20) that

$$P(I_n = 0) \leq \exp\left(-\frac{(4b_{1,2} + 4\omega(1/\sqrt{n}, \mathcal{G}))^2}{2v_{1,2}}\right)$$
$$+ \exp\left(-\frac{(4b_{2,1} + 4\omega(1/\sqrt{n}, \mathcal{G}))^2}{2v_{2,1}}\right).$$

First note that if $\omega^2(\frac{1}{\sqrt{n}}, \mathcal{F}_1) \geq \frac{1}{\ln \gamma_{1,2}} \omega^2(\sqrt{\frac{\ln \gamma_{1,2}}{n}}, \mathcal{F}_1, \mathcal{F}_2)$, then since $e^{-2x} < \frac{1}{2} x^{-2}$ for $x > 0$, it follows that

$$\exp\left(-\frac{(4b_{1,2} + 4\omega(1/\sqrt{n}, \mathcal{G}))^2}{2v_{1,2}}\right) \leq \exp\left(-2\frac{\omega^2(1/\sqrt{n}, \mathcal{G})}{\omega^2(1/\sqrt{n}, \mathcal{F}_1)}\right)$$
$$\leq \frac{1}{2} \frac{\omega^4(1/\sqrt{n}, \mathcal{F}_1)}{\omega^4(1/\sqrt{n}, \mathcal{G})}. \tag{26}$$

On the other hand, if $\omega^2(\frac{1}{\sqrt{n}}, \mathcal{F}_1) < \frac{1}{\ln \gamma_{1,2}} \omega^2(\sqrt{\frac{\ln \gamma_{1,2}}{n}}, \mathcal{F}_1, \mathcal{F}_2)$, then

$$\exp\left(-\frac{(4b_{1,2} + 4\omega(1/\sqrt{n}, \mathcal{G}))^2}{2v_{1,2}}\right)$$
$$\leq \exp\left(-\frac{16\omega^2(\sqrt{\ln \gamma_{1,2}/n}, \mathcal{F}_1, \mathcal{F}_2) + 16\omega^2(1/\sqrt{n}, \mathcal{G})}{(8/\ln \gamma_{1,2})\omega^2(\sqrt{\ln \gamma_{1,2}/n}, \mathcal{F}_1, \mathcal{F}_2)}\right)$$
$$\leq \exp\left(-\left(4\ln \gamma_{1,2} + 2\frac{\omega^2(1/\sqrt{n}, \mathcal{G})}{\omega^2(1/\sqrt{n}, \mathcal{F}_1, \mathcal{F}_2)}\right)\right) \tag{27}$$

$$\leq \exp\left(-\left(4\ln\gamma_{1,2} + 2\frac{\omega^2(1/\sqrt{n},\mathcal{G})}{\omega^2(1/\sqrt{n},\mathcal{F}_1,\mathcal{G})}\right)\right)$$

$$\leq \frac{\omega^4(1/\sqrt{n},\mathcal{F}_1)}{\omega^4(1/\sqrt{n},\mathcal{F}_1,\mathcal{G})} \cdot \frac{1}{2}\frac{\omega^4(1/\sqrt{n},\mathcal{F}_1,\mathcal{G})}{\omega^4(1/\sqrt{n},\mathcal{G})} = \frac{1}{2}\frac{\omega^4(1/\sqrt{n},\mathcal{F}_1)}{\omega^4(1/\sqrt{n},\mathcal{G})}.$$

Hence combining (26) and (27) yields $\exp(-\frac{(4b_{1,2}+4\omega(1/\sqrt{n},\mathcal{G}))^2}{2v_{1,2}}) \leq \frac{1}{2}\frac{\omega^4(1/\sqrt{n},\mathcal{F}_1)}{\omega^4(1/\sqrt{n},\mathcal{G})}$.
A similar argument shows $\exp(-\frac{(4b_{2,1}+4\omega(1/\sqrt{n},\mathcal{G}))^2}{2v_{2,1}}) \leq \frac{1}{2}\frac{\omega^4(1/\sqrt{n},\mathcal{F}_1)}{\omega^4(1/\sqrt{n},\mathcal{G})}$ and (25) follows. □

The test also has a large probability of rejecting $\mathcal{F}_1$ when $f \in \mathcal{F}_2$ and the bias of $\hat{T}_1$ is large since in such a case either $E(\hat{T}_1 - \hat{T}_{2,1})$ is large or $E(\hat{T}_{1,2} - \hat{T}_1)$ is large. The following lemma gives a useful upper bound on the probability of using $\hat{T}_1$ in this case.

LEMMA 2. *If $f \in \mathcal{F}_2$ and $|E\hat{T}_1 - Tf| \geq \lambda(b_{1,2} + b_{2,1} + \omega(\frac{1}{\sqrt{n}},\mathcal{G}))$ for some $\lambda > 6$, then*

$$(28) \qquad P(I_n = 1) \leq e^{-(\lambda-6)^2/4}.$$

PROOF. We shall only give the proof when $E\hat{T}_1 - Tf \geq \lambda(b_{1,2} + b_{2,1} + \omega(\frac{1}{\sqrt{n}},\mathcal{G}))$, as the case when $E\hat{T}_1 - Tf \leq -\lambda(b_{1,2} + b_{2,1} + \omega(\frac{1}{\sqrt{n}},\mathcal{G}))$ can be handled similarly. Let $f \in \mathcal{F}_2$. Then $P(I_n = 1) \leq P(\hat{T}_1 - \hat{T}_{2,1} \leq 5b_{2,1} + 4\omega(\frac{1}{\sqrt{n}},\mathcal{G}))$. Note that

$$E\left(\hat{T}_1 - \hat{T}_{2,1} - 5b_{2,1} - 4\omega\left(\frac{1}{\sqrt{n}},\mathcal{G}\right)\right)$$

$$= E(\hat{T}_1 - Tf) - E(\hat{T}_{2,1} - Tf) - 5b_{2,1} - 4\omega\left(\frac{1}{\sqrt{n}},\mathcal{G}\right)$$

$$\geq \lambda b_{2,1} + \lambda\omega\left(\frac{1}{\sqrt{n}},\mathcal{G}\right) - \frac{1}{2}\omega\left(\sqrt{\frac{\ln\gamma_{2,1}}{n}},\mathcal{F}_2,\mathcal{F}_1\right) - 5b_{2,1} - 4\omega\left(\frac{1}{\sqrt{n}},\mathcal{G}\right)$$

$$\geq (\lambda - 6)\left(\omega\left(\sqrt{\frac{\ln\gamma_{2,1}}{n}},\mathcal{F}_2,\mathcal{F}_1\right) + \omega\left(\frac{1}{\sqrt{n}},\mathcal{G}\right)\right).$$

Now $\text{Var}(\hat{T}_1 - \hat{T}_{2,1}) \leq v_{2,1} = 2(\omega^2(\frac{1}{\sqrt{n}},\mathcal{F}_1) + \frac{1}{\ln\gamma_{2,1}}\omega^2(\sqrt{\frac{\ln\gamma_{2,1}}{n}},\mathcal{F}_2,\mathcal{F}_1))$ yields

$$P(I_n = 1) \leq \exp\left(-\frac{(\lambda-6)^2}{2}\frac{(\omega(\sqrt{\ln\gamma_{2,1}/n},\mathcal{F}_2,\mathcal{F}_1) + \omega(1/\sqrt{n},\mathcal{G}))^2}{v_{2,1}}\right)$$

$$\leq \exp\left(-\frac{(\lambda-6)^2}{4}\right).$$

□

The proof of Theorem 3 now follows from Lemma 1 for (22) and from Lemma 2 for (23).

PROOF OF THEOREM 3. The minimax rate optimality of $\hat{T}$ over $\mathcal{F}_1$ follows directly from Lemma 1 and the fact that $\hat{T}_2^*$ satisfies (11):

$$\sup_{f \in \mathcal{F}_1} E(\hat{T} - Tf)^2 \leq \sup_{f \in \mathcal{F}_1} E(\hat{T}_1 - Tf)^2$$
$$+ \sup_{f \in \mathcal{F}_1} (E|\hat{T}_2^* - Tf|^4)^{1/2} \cdot (P(I_n = 0))^{1/2}$$
$$\leq \omega^2\left(\frac{1}{\sqrt{n}}, \mathcal{F}_1\right) + C\omega^2\left(\frac{1}{\sqrt{n}}, \mathcal{G}\right) \frac{\omega^2(1/\sqrt{n}, \mathcal{F}_1)}{\omega^2(1/\sqrt{n}, \mathcal{G})}$$
$$= C\omega^2\left(\frac{1}{\sqrt{n}}, \mathcal{F}_1\right)$$

and thus (22) holds. The proof of (23) is broken into two parts. If $f \in \mathcal{F}_2$ and $|E\hat{T}_1 - Tf| \leq 6(b_{1,2} + b_{2,1} + \omega(\frac{1}{\sqrt{n}}, \mathcal{G}))$, then

$$E(\hat{T} - Tf)^2 \leq E(\hat{T}_1 - Tf)^2 + E(\hat{T}_2^* - Tf)^2$$
$$\leq \omega^2\left(\frac{1}{\sqrt{n}}, \mathcal{F}_1\right) + 36\left(b_{1,2} + b_{2,1} + \omega\left(\frac{1}{\sqrt{n}}, \mathcal{G}\right)\right)^2$$
$$+ C\omega^2\left(\frac{1}{\sqrt{n}}, \mathcal{G}\right)$$
(29)
$$\leq C\omega^2\left(\sqrt{\frac{\ln \gamma_{1,2}}{n}}, \mathcal{F}_1, \mathcal{F}_2\right) + C\omega^2\left(\sqrt{\frac{\ln \gamma_{2,1}}{n}}, \mathcal{F}_2, \mathcal{F}_1\right)$$
$$+ C\omega^2\left(\frac{1}{\sqrt{n}}, \mathcal{G}\right)$$
$$\leq C\left(\omega_+^2\left(\sqrt{\frac{\ln \gamma_+}{n}}, \mathcal{F}_1, \mathcal{F}_2\right) + \omega^2\left(\frac{1}{\sqrt{n}}, \mathcal{F}_2\right)\right),$$

where $C$ is a constant not depending on $f$, and hence in this case (23) holds.

Now note that if $X$ has a normal distribution with mean $\mu$ and variance $\sigma^2$, then

(30) $\qquad (EX^4)^{1/2} \leq 3(\mu^2 + \sigma^2).$

Hence if $f \in \mathcal{F}_2$ and $|E\hat{T}_1 - Tf| \geq \lambda(b_{1,2} + b_{2,1} + \omega(\frac{1}{\sqrt{n}}, \mathcal{G}))$ for some $\lambda > 6$, it then follows from Lemma 2 and inequality (30) that

$$E(\hat{T} - Tf)^2 \leq (E|\hat{T}_1 - Tf|^4)^{1/2} \cdot (P(I_n = 1))^{1/2} + E(\hat{T}_2^* - Tf)^2$$

$$
\begin{aligned}
&\leq \left(3\operatorname{Var}(\hat{T}_1) + 3\lambda^2\left(b_{1,2} + b_{2,1} + \omega\left(\frac{1}{\sqrt{n}}, \mathcal{G}\right)\right)^2\right) \cdot e^{-(\lambda-6)^2/8} \\
&\quad + C\omega^2\left(\frac{1}{\sqrt{n}}, \mathcal{G}\right) \\
&\leq C\omega^2\left(\sqrt{\frac{\ln \gamma_{1,2}}{n}}, \mathcal{F}_1, \mathcal{F}_2\right) + C\omega^2\left(\sqrt{\frac{\ln \gamma_{2,1}}{n}}, \mathcal{F}_2, \mathcal{F}_1\right) \\
&\quad + C\omega^2\left(\frac{1}{\sqrt{n}}, \mathcal{G}\right) \\
&\leq C\left(\omega_+^2\left(\sqrt{\frac{\ln \gamma_+}{n}}, \mathcal{F}_1, \mathcal{F}_2\right) + \omega^2\left(\frac{1}{\sqrt{n}}, \mathcal{F}_2\right)\right),
\end{aligned}
\tag{31}
$$

where the constant $C$ does not depend on $f$, and so (23) also holds in this case and the theorem follows. $\square$

**3. Examples.** Section 2 develops the general theory of optimally adaptive estimation over two convex parameter spaces. The results can be usefully explained in an alternative way. Let $\mathcal{F}_1$ and $\mathcal{F}_2$ be two convex parameter spaces with nonempty intersection and $\omega(\varepsilon, \mathcal{F}_1) \leq \omega(\varepsilon, \mathcal{F}_2)$ for $0 \leq \varepsilon \leq \varepsilon_0$. Let $\mathcal{T}_{n,c}(\mathcal{F}_1)$ be the collection of estimators which satisfy

$$\mathcal{T}_{n,c}(\mathcal{F}_1) = \left\{\hat{T}: \sup_{f \in \mathcal{F}_1} E_f(\hat{T} - Tf)^2 \leq c^2\omega^2\left(\frac{1}{\sqrt{n}}, \mathcal{F}_1\right)\right\}$$

and let

$$R_{n,c}(\mathcal{F}_1, \mathcal{F}_2) = \inf_{\hat{T} \in \mathcal{T}_{n,c}(\mathcal{F}_1)} \sup_{f \in \mathcal{F}_2} E_f(\hat{T} - Tf)^2. \tag{32}$$

The quantity $R_{n,c}(\mathcal{F}_1, \mathcal{F}_2)$ gives the optimal performance over $\mathcal{F}_2$ for minimax rate optimal estimators over $\mathcal{F}_1$. Theorems 1 and 3 taken together quantify $R_{n,c}(\mathcal{F}_1, \mathcal{F}_2)$ in terms of the between-class modulus of continuity as

$$R_{n,c}(\mathcal{F}_1, \mathcal{F}_2) \asymp \omega_+^2\left(\sqrt{\frac{\ln \gamma_+}{n}}, \mathcal{F}_1, \mathcal{F}_2\right) + \omega^2\left(\frac{1}{\sqrt{n}}, \mathcal{F}_2\right), \tag{33}$$

where $a_n \asymp b_n$ means that $a_n/b_n$ is bounded away from 0 and $\infty$ as $n \to \infty$.

In most common cases when estimating a linear functional over convex parameter spaces the moduli are Hölderian,

$$\omega_+(\varepsilon, \mathcal{F}_i, \mathcal{F}_j) = C_{i,j}\varepsilon^{q(\mathcal{F}_i, \mathcal{F}_j)}(1 + o(1)), \tag{34}$$

where we shall write $q(\mathcal{F}_i)$ for $q(\mathcal{F}_i, \mathcal{F}_i)$. In such cases especially clear and precise statements can be made which are direct consequences of Theorems 1 and 3.

COROLLARY 1. *Let $\mathcal{F}_1$ and $\mathcal{F}_2$ be convex parameter spaces with $\mathcal{F}_1 \cap \mathcal{F}_2 \neq \varnothing$ and let $T$ be a linear functional. Suppose $\omega_+(\varepsilon, \mathcal{F}_i, \mathcal{F}_j)$ are Hölderian with exponent $q(\mathcal{F}_i, \mathcal{F}_j)$ for $i, j = 1, 2$. If $q(\mathcal{F}_1, \mathcal{F}_2) = q(\mathcal{F}_2) < q(\mathcal{F}_1)$ or $q(\mathcal{F}_1, \mathcal{F}_2) < q(\mathcal{F}_2) \leq q(\mathcal{F}_1)$, then*

$$(35) \quad C_1 \left( \frac{\log n}{n} \right)^{q(\mathcal{F}_1, \mathcal{F}_2)} \leq R_{n,c}(\mathcal{F}_1, \mathcal{F}_2) \leq C_2 \left( \frac{\log n}{n} \right)^{q(\mathcal{F}_1, \mathcal{F}_2)},$$

*where $0 < C_1 \leq C_2$ are constants and $R_{n,c}(\mathcal{F}_1, \mathcal{F}_2)$ is defined above in* (32).

Corollary 1 can then be used to classify the problem of adaptation over convex parameter spaces into three cases:

• *Case* 1. $q(\mathcal{F}_1, \mathcal{F}_2) = \min(q(\mathcal{F}_1), q(\mathcal{F}_2)) < \max(q(\mathcal{F}_1), q(\mathcal{F}_2))$. This is the "regular case" which holds for many linear functionals and common function classes of interest. In this case, one must lose a logarithmic factor as the minimum cost for adaptation. A common example of such a case is estimating a function or a derivative at a point, that is, $Tf = f^{(s)}(t_0)$ for some $s \geq 0$ when the parameter spaces are assumed to be Lipschitz. See Example 2 below and [3, 14, 17].

Besides the regular case, there are two extreme cases.

• *Case* 2. $q(\mathcal{F}_1, \mathcal{F}_2) > \min(q(\mathcal{F}_1), q(\mathcal{F}_2))$ or $q(\mathcal{F}_1, \mathcal{F}_2) = q(\mathcal{F}_1) = q(\mathcal{F}_2)$. This is a case which is not covered in Corollary 1. Results given in Section 2 show that in this case adaptation for free is always possible. That is, one can attain the optimal rate of convergence over $\mathcal{F}_1$ and $\mathcal{F}_2$ simultaneously. An example of this case is estimating a function at a point over two monotone Lipschitz classes. See Examples 1 and 3 below and [20].

• *Case* 3. $q(\mathcal{F}_1, \mathcal{F}_2) < \min(q(\mathcal{F}_1), q(\mathcal{F}_2))$. In this case the cost of adaptation is significant, much more than the usual logarithmic penalty in the regular case. If $f$ is known to be in $\mathcal{F}_1$, one can attain the rate of $n^{q(\mathcal{F}_1)}$; and if one knows that $f$ is in $\mathcal{F}_2$, the rate of convergence $n^{q(\mathcal{F}_2)}$ can be achieved. Without the information, however, one can only achieve the rate of $(n/\log n)^{q(\mathcal{F}_1, \mathcal{F}_2)}$ at best. So the cost of adaptation is a power of $n$ rather than the logarithmic factor as in the regular case. See Example 2 below.

Note that if the the parameter spaces $\mathcal{F}_1$ and $\mathcal{F}_2$ are nested, then only Cases 1 and 2 are possible and Case 3 does not arise.

We now consider a few examples below to illustrate the three different cases. Examples 1 and 3 cover Case 2 in which full adaptation is possible. Example 2 covers both Case 1 and Case 3 with different choices of parameters. In each of these examples we need to calculate the between-class modulus of continuity. The basic idea behind these calculations is contained in [9] and consists of finding extremal functions. See [4] for the details of these calculations.

EXAMPLE 1. In this example we shall have $0 < q(\mathcal{F}_2) < q(\mathcal{F}_1, \mathcal{F}_2) = q(\mathcal{F}_1) < 1$ and $\omega(\varepsilon, \mathcal{F}_1, \mathcal{F}_2) = \omega(\varepsilon, \mathcal{F}_2, \mathcal{F}_1)$. In this case full mean squared error adaptation is possible.

For $0 < \alpha \leq 1$, let

$$F(\alpha, M) = \{f : [-\tfrac{1}{2}, \tfrac{1}{2}] \to \mathbb{R} : |f(x) - f(y)| \leq M|x - y|^\alpha\}. \tag{36}$$

Let $D$ be the set of all decreasing functions and let $F_D(\alpha, M) = F(\alpha, M) \cap D$ be the set of decreasing functions which are also members of $F(\alpha, M)$. Let $Tf = f(0)$ and assume that $0 < \alpha_2 < \alpha_1 \leq 1$. Let $\mathcal{F}_1 = F_D(\alpha_1, M_1)$ and $\mathcal{F}_2 = F_D(\alpha_2, M_2)$. Then for these parameter spaces and the linear functional $Tf = f(0)$ it follows from calculations in [4] that, as $\varepsilon \to 0$,

$$\omega^2(\varepsilon, \mathcal{F}_1, \mathcal{F}_2) = \omega^2(\varepsilon, \mathcal{F}_2, \mathcal{F}_1)$$
$$= (2\alpha_1 + 1)^{\alpha_1/(2\alpha_1+1)} M_1^{1/(2\alpha_1+1)} \varepsilon^{2\alpha_1/(2\alpha_1+1)} (1 + o(1)), \tag{37}$$

$$\omega^2(\varepsilon, \mathcal{F}_1) = (\alpha_1 + \tfrac{1}{2})^{\alpha_1/(2\alpha_1+1)} M_1^{1/(2\alpha_1+1)} \varepsilon^{2\alpha_1/(2\alpha_1+1)} (1 + o(1)), \tag{38}$$

$$\omega^2(\varepsilon, \mathcal{F}_2) = (\alpha_2 + \tfrac{1}{2})^{\alpha_2/(2\alpha_2+1)} M_2^{1/(2\alpha_2+1)} \varepsilon^{2\alpha_2/(2\alpha_2+1)} (1 + o(1)). \tag{39}$$

In this case $q(\mathcal{F}_1, \mathcal{F}_2) = \max(q(\mathcal{F}_1), q(\mathcal{F}_2)) > \min(q(\mathcal{F}_1), q(\mathcal{F}_2))$ and hence adaptation for free can be achieved.

EXAMPLE 2. This example shows that sometimes we must lose more than a logarithmic factor when we try to adapt. Let

$$F^R(\alpha, M) = \{f : [-\tfrac{1}{2}, \tfrac{1}{2}] \to \mathbb{R} : |f^{(s)}(x) - f^{(s)}(y)| \leq M|x - y|^{\alpha-s}$$
$$0 \leq x \leq y \leq \tfrac{1}{2}\},$$

where $s$ is the largest integer less than $\alpha$. Similarly, let

$$F^L(\alpha, M) = \{f : [-\tfrac{1}{2}, \tfrac{1}{2}] \to \mathbb{R} : |f^{(s)}(x) - f^{(s)}(y)| \leq M|x - y|^{\alpha-s}$$
$$-\tfrac{1}{2} \leq x \leq y \leq 0\}.$$

Finally let $F(\alpha_1, M_1, \alpha_2, M_2) = F^L(\alpha_1, M_1) \cap F^R(\alpha_2, M_2)$.

Note that for the linear functional $Tf = f(0)$ and the (ordered) parameter spaces $\mathcal{F}_1 = F(\alpha_1, M_1, \alpha_2, M_2)$ and $\mathcal{F}_2 = F(\beta_1, N_1, \beta_2, N_2)$ it follows from the calculations given in [4] that

$$\omega^2(\varepsilon, \mathcal{F}_1) = C(\alpha_1, M_1, \alpha_2, M_2) \varepsilon^{2\delta/(2\delta+1)} (1 + o(1)), \tag{40}$$

$$\omega^2(\varepsilon, \mathcal{F}_2) = C(\beta_1, N_1, \beta_2, N_2) \varepsilon^{2\rho/(2\rho+1)} (1 + o(1)), \tag{41}$$

where $\delta = \max(\alpha_1, \alpha_2)$ and $\rho = \max(\beta_1, \beta_2)$.

Now let $0 < \alpha_2 \leq \alpha_1 \leq 1$ and $0 < \beta_1 \leq \beta_2 \leq 1$. Then $q(\mathcal{F}_1) = \frac{2\alpha_1}{2\alpha_1+1}$ and $q(\mathcal{F}_2) = \frac{2\beta_2}{2\beta_2+1}$. The between-class modulus satisfies

$$\omega^2(\varepsilon, \mathcal{F}_1, \mathcal{F}_2) = C(M_1, \alpha_1, M_2, \alpha_2, N_1, \beta_1, N_2, \beta_2) \varepsilon^{2\gamma/(2\gamma+1)} (1 + o(1)), \tag{42}$$

where $\gamma = \max(\min(\alpha_1, \beta_1), \min(\alpha_2, \beta_2))$.

Two interesting cases arise, depending on the relationship among $\alpha_1, \alpha_2, \beta_1$ and $\beta_2$.

- $\beta_2 > \beta_1 \geq \alpha_1 > \alpha_2$. Then the quantity $\gamma$ in (42) is $\gamma = \alpha_1$ and so $q(\mathcal{F}_1, \mathcal{F}_2) = \frac{2\alpha_1}{2\alpha_1+1}$. Hence in this case

$$q(\mathcal{F}_1, \mathcal{F}_2) = \min(q(\mathcal{F}_1), q(\mathcal{F}_2)) < \max(q(\mathcal{F}_1), q(\mathcal{F}_2)).$$

This is a case where a logarithmic penalty term must be paid for adaptation.

- $\alpha_1 \geq \beta_2 > \beta_1 \geq \alpha_2$. In this case, the quantity $\gamma$ in (42) is $\gamma = \beta_1$ and hence $q(\mathcal{F}_1, \mathcal{F}_2) = \frac{2\beta_1}{2\beta_1+1}$. Therefore in this case

$$q(\mathcal{F}_1, \mathcal{F}_2) < \min(q(\mathcal{F}_1), q(\mathcal{F}_2)).$$

Consequently the cost of adaption between $\mathcal{F}_1$ and $\mathcal{F}_2$ is much more than a logarithmic penalty. The maximum risk over the two spaces is of the order $n^{-2\beta_1/(2\beta_1+1)}$.

A particularly interesting case is when $\alpha_1 = \beta_2 > \beta_1 \geq \alpha_2$. In this case, the minimax rates of convergence over $\mathcal{F}_1$ and $\mathcal{F}_2$ are the same, both equal to $n^{-2\beta_2/(2\beta_2+1)}$. Yet it is impossible to achieve this optimal rate adaptively over the two parameter spaces; in fact the cost of adaption in this case is substantial.

EXAMPLE 3. This will give an example where $0 < q(\mathcal{F}_1) < q(\mathcal{F}_1, \mathcal{F}_2) < q(\mathcal{F}_2) < 1$. It will also yield an example where $\omega(\varepsilon, \mathcal{F}_1, \mathcal{F}_2) \neq \omega(\varepsilon, \mathcal{F}_2, \mathcal{F}_1)$. In this case full mean squared error adaptation can be achieved. Let $Tf = f(0)$. Now let

$$F_D(\alpha_1, M_1, \alpha_2, M_2) = F(\alpha_1, M_1, \alpha_2, M_2) \cap D,$$

where $F(\alpha_1, M_1, \alpha_2, M_2)$ is defined as in Example 2.

Let $\beta_1 > \beta_2 > \alpha_1 > \alpha_2$. Calculations in [4] yield for the (ordered) parameter spaces $\mathcal{F}_1 = F_D(\alpha_1, M_1, \alpha_2, M_2)$ and $\mathcal{F}_2 = F_D(\beta_1, N_1, \beta_2, N_2)$,

$$\omega^2(\varepsilon, \mathcal{F}_1) = C\varepsilon^{2\alpha_1/(2\alpha_1+1)}(1 + o(1)),$$

(43)

$$\omega^2(\varepsilon, \mathcal{F}_2) = C\varepsilon^{2\beta_1/(2\beta_1+1)}(1 + o(1)),$$

$$\omega^2(\varepsilon, \mathcal{F}_1, \mathcal{F}_2) = C\varepsilon^{2\beta_2/(2\beta_2+1)}(1 + o(1)),$$

(44)

$$\omega^2(\varepsilon, \mathcal{F}_2, \mathcal{F}_1) = C\varepsilon^{2\beta_1/(2\beta_1+1)}(1 + o(1)).$$

Hence this is an example where $\omega(\varepsilon, \mathcal{F}_1, \mathcal{F}_2) \neq \omega(\varepsilon, \mathcal{F}_2, \mathcal{F}_1)(1 + o(1))$. Note that $\beta_1 > \beta_2$. It then follows from (44) that $q(\mathcal{F}_1, \mathcal{F}_2) = \frac{2\beta_2}{2\beta_2+1}$. Hence this is an example where

$$0 < q(\mathcal{F}_1) < q(\mathcal{F}_1, \mathcal{F}_2) < q(\mathcal{F}_2) < 1.$$

In particular, $q(\mathcal{F}_1, \mathcal{F}_2) > \min(q(\mathcal{F}_1), q(\mathcal{F}_2))$, so it is also an example where full mean squared error adaptation is possible.

**4. Adaptation over many parameter spaces.** Section 2 gives a complete treatment of adaptation over two convex parameter spaces. It is shown that the between-class modulus determines the cost of adaptation and the ordered modulus can be used for the construction of optimally adaptive procedures. This theory of adaptation over two parameter spaces is in turn a fundamental building block for adaptation over richer collections of parameter spaces. We first extend the theory to adaptation over collections of finitely many parameter spaces. Section 5 further generalizes the theory to collections of infinitely many parameter spaces.

The basic idea for the construction of adaptive estimators builds on that given for two parameter spaces. In particular the adaptive estimator is based on the construction of tests between pairs of parameter spaces. The resulting estimator is optimally adaptive in the sense defined in Section 2: it attains the lower bound on the cost of adaptation over finitely many convex parameter spaces which satisfy certain regularity conditions on the moduli. We shall begin by assuming that the parameter spaces are nested, in which case these conditions are always satisfied.

4.1. *Adaptation over nested parameter spaces.* Let $\mathcal{F}_1 \subset \mathcal{F}_2 \subset \cdots \subset \mathcal{F}_k$ be closed convex parameter spaces and for convenience set $\mathcal{F}_0 = \varnothing$. In this context the goal of adaptation is most easily described sequentially. First, the estimator should attain the exact minimax rate of convergence over $\mathcal{F}_1$. Given the performance over $\mathcal{F}_1$, the estimator should attain the lower bound as given in Theorem 1 over $\mathcal{F}_2$. Moreover, for $i \geq 3$ the estimator should attain the lower bound given its performance over $\mathcal{F}_1, \mathcal{F}_2, \ldots, \mathcal{F}_{i-1}$.

We shall introduce some notation before explaining the lower bounds in detail. For $i \neq j$, define the quantity $\gamma_{i,j} > 0$ as follows. If $i \wedge j = \min(i,j) = 1$, let

$$\gamma_{i,j}^2 = \max\left(e, \frac{\omega^2(1/\sqrt{n}, \mathcal{F}_i, \mathcal{F}_j)}{\omega^2(1/\sqrt{n}, \mathcal{F}_1)}\right) \quad \text{and}$$

(45)

$$\gamma_{i,j,+}^2 = \max\left(e, \frac{\omega_+^2(1/\sqrt{n}, \mathcal{F}_i, \mathcal{F}_j)}{\omega^2(1/\sqrt{n}, \mathcal{F}_1)}\right).$$

If $i \wedge j \geq 2$, define $\gamma_{i,j}$ and $\gamma_{i,j,+}$ recursively by

(46)
$$\gamma_{i,j}^2 = \max\left(e, \omega^2\left(\frac{1}{\sqrt{n}}, \mathcal{F}_i, \mathcal{F}_j\right)\right.$$
$$\left. \times \left(\max_{1 \leq m \leq i \wedge j - 1}\left\{\omega_+^2\left(\sqrt{\frac{\ln \gamma_{m,i \wedge j,+}}{n}}, \mathcal{F}_m, \mathcal{F}_{i \wedge j}\right)\right\}\right)\right\}$$

$$+ \omega^2\left(\frac{1}{\sqrt{n}}, \mathcal{F}_{i \wedge j}\right)\right)^{-1}\right)$$

and

$$\gamma_{i,j,+}^2 = \max(\gamma_{i,j}, \gamma_{j,i})$$
$$= \max\left(e, \omega_+^2\left(\frac{1}{\sqrt{n}}, \mathcal{F}_i, \mathcal{F}_j\right)\right.$$

(47)
$$\times \left(\max_{1 \leq m \leq i \wedge j - 1}\left\{\omega_+^2\left(\sqrt{\frac{\ln \gamma_{m,i \wedge j,+}}{n}}, \mathcal{F}_m, \mathcal{F}_{i \wedge j}\right)\right\}\right.$$
$$\left.\left. + \omega^2\left(\frac{1}{\sqrt{n}}, \mathcal{F}_{i \wedge j}\right)\right)^{-1}\right).$$

Let $A_i(n) \geq 0$ be defined by $A_1^2(n) = \omega^2(\frac{1}{\sqrt{n}}, \mathcal{F}_1)$ and, for $2 \leq i \leq k$,

(48) $$A_i^2(n) = \max_{1 \leq m \leq i-1}\left\{\omega_+^2\left(\sqrt{\frac{\ln \gamma_{m,i,+}}{n}}, \mathcal{F}_m, \mathcal{F}_i\right)\right\} + \omega^2\left(\frac{1}{\sqrt{n}}, \mathcal{F}_i\right).$$

Suppose that $c_i > 0$ are some constants for $i = 1, \ldots k$. If $\hat{T}$ is an estimator satisfying

(49) $$\sup_{\mathcal{F}_1} E(\hat{T} - Tf)^2 \leq c_1 A_1^2(n),$$

then Theorem 1 shows that the estimator $\hat{T}$ must satisfy a lower bound over $\mathcal{F}_2$,

(50) $$\sup_{\mathcal{F}_2} E(\hat{T} - Tf)^2 \geq d_2 A_2^2(n),$$

where $d_2 > 0$ is a constant. More generally for $2 \leq j \leq k$, if an estimator $\hat{T}$ satisfies

(51) $$\sup_{\mathcal{F}_i} E(\hat{T} - Tf)^2 \leq c_i A_i^2(n) \qquad \text{for } i = 1, \ldots, j-1,$$

then Theorem 2 shows that the estimator $\hat{T}$ must satisfy a lower bound over $\mathcal{F}_j$,

(52) $$\sup_{\mathcal{F}_j} E(\hat{T} - Tf)^2 \geq d_j A_j^2(n)$$

for some constant $d_j > 0$. It is thus natural to seek an estimator which attains (51) for all $1 \leq i \leq k$ for some constants $c_i > 0$. In light of the lower bound (52), such an estimator can also be termed optimally adaptive.

We now turn to the construction of such adaptive estimators. As in Section 2.2, let $\hat{T}_i$ be linear estimators satisfying $\sup_{f \in \mathcal{F}_i} E(\hat{T}_i - Tf)^2 \leq \omega^2(\frac{1}{\sqrt{n}}, \mathcal{F}_i)$. The procedure, which will be defined precisely below, can be described sequentially as follows.

1. Test between $\mathcal{F}_1$ and $\mathcal{F}_i$ for all $2 \leq i \leq k$.
2. If all the tests are in favor of $\mathcal{F}_1$, use $\hat{T}_1$ as the estimate of $Tf$.
3. Otherwise, delete $\mathcal{F}_1$ and repeat steps 1 and 2.

The performance of this procedure depends critically on the properties of the tests between pairs of parameter spaces. The tests are developed in a similar but somewhat more involved way than those in Section 2.

For $i \neq j$ let $\hat{T}_{i,j}$ be the estimator satisfying (12)–(14) with $\mathcal{F} = \mathcal{F}_i$, $\mathcal{H} = \mathcal{F}_j$ and $V = \sigma_{i,j}^2$, where $\sigma_{i,j}^2 = \frac{1}{\ln \gamma_{i,j}} \omega^2(\sqrt{\frac{\ln \gamma_{i,j}}{n}}, \mathcal{F}_i, \mathcal{F}_j)$. Then note that as in Section 2.2, if $f \in \mathcal{F}_i$,

$$(53) \quad E(\hat{T}_i - \hat{T}_{i,j}) \geq -\omega\left(\frac{1}{\sqrt{n}}, \mathcal{F}_i\right) - \omega\left(\sqrt{\frac{\ln \gamma_{i,j}}{n}}, \mathcal{F}_i, \mathcal{F}_j\right) = -b_{i,j},$$

$$(54) \quad E(\hat{T}_i - \hat{T}_{j,i}) \leq \omega\left(\frac{1}{\sqrt{n}}, \mathcal{F}_i\right) + \omega\left(\sqrt{\frac{\ln \gamma_{i,j}}{n}}, \mathcal{F}_j, \mathcal{F}_i\right) = b_{j,i},$$

$$(55) \quad \mathrm{Var}(\hat{T}_i - \hat{T}_{i,j}) \leq 2\left(\omega^2\left(\frac{1}{\sqrt{n}}, \mathcal{F}_i\right) + \frac{1}{\ln \gamma_{i,j}} \omega^2\left(\sqrt{\frac{\ln \gamma_{i,j}}{n}}, \mathcal{F}_i, \mathcal{F}_j\right)\right) = v_{i,j},$$

$$(56) \quad \mathrm{Var}(\hat{T}_i - \hat{T}_{j,i}) \leq 2\left(\omega^2\left(\frac{1}{\sqrt{n}}, \mathcal{F}_i\right) + \frac{1}{\ln \gamma_{i,j}} \omega^2\left(\sqrt{\frac{\ln \gamma_{i,j}}{n}}, \mathcal{F}_j, \mathcal{F}_i\right)\right) = v_{j,i}.$$

For $i < j$ the test between $\mathcal{F}_i$ and $\mathcal{F}_j$ is given by

$$(57) \quad I_{i,j} = 1(\hat{T}_{i,j} - (4(2k)^{1/2} + 1)b_{i,j} - 4k^{1/2} A_j(n) \leq \hat{T}_i$$
$$\leq \hat{T}_{j,i} + (4(2k)^{1/2} + 1)b_{j,i} + 4k^{1/2} A_j(n)).$$

The test is in favor of $\mathcal{F}_i$ if $I_{i,j} = 1$. Our adaptive estimation procedure is defined in terms of the tests $I_{i,j}$ and the minimax rate optimal estimator over $\mathcal{F}_i$, $\hat{T}_i$. The procedure is defined sequentially from "inside–out." It first tests if $f \in \mathcal{F}_1$ by checking whether $\prod_{j \geq 2} I_{1,j} = 1$, which means that all the tests $I_{1,j}$ are in favor of $\mathcal{F}_1$. In this case $\hat{T}_1$ is used. Otherwise $\mathcal{F}_1$ is deleted and the procedure iterates. More formally, the estimator $\hat{T}^*$ can be written as

$$(58) \quad \hat{T}^* = \sum_{i=1}^{k} \left(1 - \prod_{m<i} \prod_{j>m} I_{m,j}\right)\left(\prod_{j>i} I_{i,j}\right) \hat{T}_i.$$

The following theorem shows that this procedure is optimally adaptive in the sense that it attains (51) for all $1 \leq i \leq k$.

THEOREM 4. *Let $\mathcal{F}_1 \subset \mathcal{F}_2 \subset \cdots \subset \mathcal{F}_k$ be closed convex parameter spaces. Then the estimator $\hat{T}^*$ defined in (58) is optimally adaptive over $\mathcal{F}_i$ for $i = 1, \ldots, k$. More specifically,*

$$\sup_{f \in \mathcal{F}_1} E(\hat{T}^* - Tf)^2 \leq C\omega^2\left(\frac{1}{\sqrt{n}}, \mathcal{F}_1\right), \tag{59}$$

*and for $2 \leq i \leq k$,*

$$\sup_{f \in \mathcal{F}_i} E(\hat{T}^* - Tf)^2 \tag{60}$$
$$\leq C\left(\max_{1 \leq m \leq i-1}\left\{\omega_+^2\left(\sqrt{\frac{\ln \gamma_{m,i,+}}{n}}, \mathcal{F}_m, \mathcal{F}_i\right)\right\} + \omega^2\left(\frac{1}{\sqrt{n}}, \mathcal{F}_i\right)\right).$$

The basic ideas for the proof of Theorem 4 are similar to those of Theorem 3, but the calculations involved are more complicated. There are two main concerns which need to be addressed. For $f \in \mathcal{F}_i \setminus \mathcal{F}_{i-1}$ one concern is that the test stops too late and uses $\hat{T}_j$ for some $j > i$. Lemma 3 below shows that this probability is small. The other concern is that the test stops early and uses $\hat{T}_j$ for $j < i$. This is only a problem when the bias of $\hat{T}_j$ is large. We shall show that if that is indeed the case, then the chance of using such a $\hat{T}_j$ is small. The specific bound is given in Lemma 4.

LEMMA 3. *If $f \in \mathcal{F}_i$, then for $j > i$,*

$$P(\hat{T}^* = \hat{T}_j) \leq k\frac{A_i^4(n)}{A_j^4(n)}, \tag{61}$$

*where $A_i(n)$ is defined as in (48). In particular, (61) holds for $f \in \mathcal{F}_i \setminus \mathcal{F}_{i-1}$.*

PROOF. It follows from (53)–(56) that for $f \in \mathcal{F}_i$ and $i < m \leq k$

$$P(I_{i,m} = 0) \leq P(\hat{T}_i - \hat{T}_{i,m} \leq -(4(2k)^{1/2} + 1)b_{i,m} - 4k^{1/2}A_m(n))$$
$$+ P(\hat{T}_i - \hat{T}_{m,i} \geq (4(2k)^{1/2} + 1)b_{m,i} + 4k^{1/2}A_m(n))$$
$$\leq \exp\left(-\frac{(4(2k)^{1/2}b_{i,m} + 4k^{1/2}A_m(n))^2}{2v_{i,m}}\right)$$
$$+ \exp\left(-\frac{(4(2k)^{1/2}b_{m,i} + 4k^{1/2}A_m(n))^2}{2v_{m,i}}\right).$$

First note that if $\omega^2(\frac{1}{\sqrt{n}}, \mathcal{F}_i) \geq \frac{1}{\ln \gamma_{i,m}} \omega^2(\sqrt{\frac{\ln \gamma_{i,m}}{n}}, \mathcal{F}_i, \mathcal{F}_m)$, then $v_{i,m} \leq 4A_i^2(n)$. Since $e^{-2kx} < \frac{1}{2}x^{-2k}$ for $x > 0$ and $k \geq 1$, it follows that

$$\exp\left(-\frac{(4(2k)^{1/2}b_{i,m} + 4k^{1/2}A_m(n))^2}{2v_{i,m}}\right) \leq \exp\left(-2k\frac{A_m^2(n)}{A_i^2(n)}\right) \leq \frac{1}{2}\frac{A_i^{4k}(n)}{A_m^{4k}(n)}.$$
(62)

On the other hand, if $\omega^2(\frac{1}{\sqrt{n}}, \mathcal{F}_i) < \frac{1}{\ln \gamma_{i,m}} \omega^2(\sqrt{\frac{\ln \gamma_{i,m}}{n}}, \mathcal{F}_i, \mathcal{F}_m)$, then

$$\exp\left(-\frac{(4(2k)^{1/2}b_{i,m} + 4k^{1/2}A_m(n))^2}{2v_{i,m}}\right)$$

$$\leq \exp\left(-\frac{32k\omega^2(\sqrt{\ln \gamma_{i,m}/n}, \mathcal{F}_i, \mathcal{F}_m) + 16kA_m^2(n)}{(8/\ln \gamma_{i,m})\omega^2(\sqrt{\ln \gamma_{i,m}/n}, \mathcal{F}_i, \mathcal{F}_m)}\right)$$

(63)
$$\leq \exp\left(-\left(4k \ln \gamma_{i,m} + 2k\frac{A_m^2(n)}{\omega^2(1/\sqrt{n}, \mathcal{F}_i, \mathcal{F}_m)}\right)\right)$$

$$\leq \frac{A_i^{4k}(n)}{\omega^{4k}(1/\sqrt{n}, \mathcal{F}_i, \mathcal{F}_m)} \cdot \frac{1}{2}\frac{\omega^{4k}(1/\sqrt{n}, \mathcal{F}_i, \mathcal{F}_m)}{A_m^{4k}(n)}$$

$$= \frac{1}{2}\frac{A_i^{4k}(n)}{A_m^{4k}(n)}.$$

Combining (62) and (63) yields $\exp(-\frac{(4(2k)^{1/2}b_{i,m} + 4k^{1/2}A_m(n))^2}{2v_{i,m}}) \leq \frac{1}{2}\frac{A_i^{4k}(n)}{A_m^{4k}(n)}$.
A similar argument yields $\exp(-\frac{(4(2k)^{1/2}b_{i,m} + 4k^{1/2}A_m(n))^2}{2v_{i,m}}) \leq \frac{1}{2}\frac{A_i^{4k}(n)}{A_m^{4k}(n)}$. Therefore,

(64)  $$P(I_{i,m} = 0) \leq \frac{A_i^{4k}(n)}{A_m^{4k}(n)} \quad \text{for } 1 \leq i < m \leq k.$$

Now note that for $j < m \leq i$, $\gamma_{j,m,+} \leq \gamma_{j,i,+}$ and consequently $\omega_+(\sqrt{\frac{\ln \gamma_{j,m,+}}{n}}, \mathcal{F}_j, \mathcal{F}_m) \leq \omega_+(\sqrt{\frac{\ln \gamma_{j,i,+}}{n}}, \mathcal{F}_j, \mathcal{F}_i)$. It then follows that $A_i(n)$ are nondecreasing in $i$ and from (64) that

(65)  $$P(I_m = 0) \leq \sum_{l=m+1}^{k} P(I_{m,l} = 0) \leq \sum_{l=m+1}^{k} \frac{A_m^{4k}(n)}{A_l^{4k}(n)} \leq k\frac{A_m^{4k}(n)}{A_{m+1}^{4k}(n)}.$$

Set $I_i = \prod_{j>i} I_{i,j}$. Then

(66)  $$P(\hat{T}^* = \hat{T}_j) \leq \min_{i \leq m \leq j-1} P(I_m = 0) \leq \left(\prod_{m=i}^{j-1} P(I_m = 0)\right)^{1/(j-i-1)}.$$

By combining (65) and (66), it follows $P(\hat{T}^* = \hat{T}_j) \leq k\frac{A_i^4(n)}{A_j^4(n)}$. □

LEMMA 4. *Suppose $f \in \mathcal{F}_i \setminus \mathcal{F}_{i-1}$ and $j < i$. If $|E\hat{T}_j - Tf| \geq \lambda(b_{j,i} + b_{i,j} + A_i(n))$ for some $\lambda > 4(2k)^{1/2} - 2$, then*

$$P(\hat{T}^* = \hat{T}_j) \leq \exp\left(-\frac{(\lambda - 4(2k)^{1/2} - 2)^2}{4}\right). \tag{67}$$

PROOF. We shall only consider the case when $E\hat{T}_j - Tf \geq \lambda(b_{j,i} + b_{i,j} + A_i(n))$ since the case of $E\hat{T}_j - Tf \leq -\lambda(b_{j,i} + b_{i,j} + A_i(n))$ can be handled similarly. Let $f \in \mathcal{F}_i \setminus \mathcal{F}_{i-1}$. Then $P(\hat{T}^* = T_j) \leq P(I_{j,i} = 1) \leq P(\hat{T}_j - \hat{T}_{i,j} \leq (4(2k)^{1/2} + 1)b_{i,j} + 4k^{1/2}A_i(n))$. Note that

$$E(\hat{T}_j - \hat{T}_{i,j} - (4(2k)^{1/2} + 1)b_{i,j} - 4k^{1/2}A_i(n))$$
$$= E(\hat{T}_j - Tf) - E(\hat{T}_{i,j} - Tf) - (4(2k)^{1/2} + 1)b_{i,j} - 4k^{1/2}A_i(n)$$
$$\geq \lambda b_{i,j} + \lambda A_i(n) - \frac{1}{2}\omega\left(\sqrt{\frac{\ln \gamma_{i,j}}{n}}, \mathcal{F}_i, \mathcal{F}_j\right)$$
$$\quad - (4(2k)^{1/2} + 1)b_{i,j} - 4k^{1/2}A_i(n)$$
$$\geq (\lambda - 4(2k)^{1/2} - 2)\left(\omega\left(\sqrt{\frac{\ln \gamma_{i,j}}{n}}, \mathcal{F}_i, \mathcal{F}_j\right) + A_i(n)\right).$$

Note that $\operatorname{Var}(\hat{T}_j - \hat{T}_{i,j}) \leq v_{i,j} = 2(\omega^2(\frac{1}{\sqrt{n}}, \mathcal{F}_j) + \frac{1}{\ln \gamma_{i,j}}\omega^2(\sqrt{\frac{\ln \gamma_{i,j}}{n}}, \mathcal{F}_i, \mathcal{F}_j))$ and hence

$$P(\hat{T}^* = T_j) \leq \exp\left(-\frac{(\lambda - 4(2k)^{1/2} - 2)^2}{2} \frac{(\omega(\sqrt{\ln \gamma_{i,j}/n}, \mathcal{F}_i, \mathcal{F}_j) + A_i(n))^2}{v_{i,j}}\right)$$
$$\leq \exp\left(-\frac{(\lambda - 4(2k)^{1/2} - 2)^2}{4}\right). \quad \square$$

We are now ready to prove Theorem 4.

PROOF OF THEOREM 4. The minimax rate optimality of $\hat{T}$ over $\mathcal{F}_1$ follows from Lemma 3:

$$\sup_{f \in \mathcal{F}_1} E(\hat{T}^* - Tf)^2 \leq \sup_{f \in \mathcal{F}_1} E(\hat{T}_1 - Tf)^2$$
$$+ \sum_{j=2}^k \sup_{f \in \mathcal{F}_1} (E|\hat{T}_j - Tf|^4)^{1/2} \cdot (P(\hat{T}^* = \hat{T}_j))^{1/2}$$
$$\leq \omega^2\left(\frac{1}{\sqrt{n}}, \mathcal{F}_1\right) + \sum_{j=2}^k \omega^2\left(\frac{1}{\sqrt{n}}, \mathcal{F}_j\right) \cdot \frac{A_1^2(n)}{A_j^2(n)}$$

$$\leq k\omega^2\left(\frac{1}{\sqrt{n}}, \mathcal{F}_1\right)$$

and hence (59) follows. The proof of (60) is somewhat more involved. Consider the case $f \in \mathcal{F}_i \setminus \mathcal{F}_{i-1}$ for some $i > 1$. Set

$$\mathcal{J}_1 = \{j < i : |E\hat{T}_j - Tf| \leq (4(2k)^{1/2} + 2)(b_{j,i} + b_{i,j} + A_i(n))\},$$

$$\mathcal{J}_2 = \{j < i : |E\hat{T}_j - Tf| > (4(2k)^{1/2} + 2)(b_{j,i} + b_{i,j} + A_i(n))\}.$$

Then for $j \in \mathcal{J}_1$, $E(\hat{T}_j - Tf)^2 \leq \omega^2(\frac{1}{\sqrt{n}}, \mathcal{F}_j) + (4(2k)^{1/2} + 2)^2(b_{j,i} + b_{i,j} + A_i(n))^2$. If $j \in \mathcal{J}_2$, then $|E\hat{T}_j - Tf| = \lambda(b_{j,i} + b_{i,j} + A_i(n))$ for some $\lambda > 4(2k)^{1/2} + 2$. Hence, by (30),

$$(E|\hat{T}_j - Tf|^4)^{1/2} \leq 3\operatorname{Var}(\hat{T}_j) + 3\lambda^2(b_{j,i} + b_{i,j} + A_i(n))^2$$

$$\leq 4\lambda^2(b_{j,i} + b_{i,j} + A_i(n))^2.$$

It then follows from Lemma 4 that $(P(\hat{T}^* = T_j))^{1/2} \leq \exp(-\frac{(\lambda - 4(2k)^{1/2} - 2)^2}{8})$. Hence, for $f \in \mathcal{F}_i \setminus \mathcal{F}_{i-1}$ with $i > 1$,

$$E(\hat{T}^* - Tf)^2 = \sum_{j=1}^{k} E\{(\hat{T}_j - Tf)^2 1(\hat{T}^* = \hat{T}_j)\}$$

$$\leq \sum_{j \in \mathcal{J}_1} E(\hat{T}_j - Tf)^2 + \sum_{j \in \mathcal{J}_2} (E|\hat{T}_j - Tf|^4)^{1/2} (P(\hat{T}^* = \hat{T}_j))^{1/2}$$

$$+ E(\hat{T}_i - Tf)^2 + \sum_{j=i+1}^{k} (E|\hat{T}_j - Tf|^4)^{1/2} (P(\hat{T}^* = \hat{T}_j))^{1/2}$$

$$\leq \sum_{j \in \mathcal{J}_1} \left\{ \omega^2\left(\frac{1}{\sqrt{n}}, \mathcal{F}_j\right) + (4(2k)^{1/2} + 2)^2(b_{j,i} + b_{i,j} + A_i(n))^2 \right\}$$

$$+ \sum_{j \in \mathcal{J}_2} 4\lambda^2(b_{j,i} + b_{i,j} + A_i(n))^2 \cdot \exp\left(-\frac{(\lambda - 4(2k)^{1/2} - 2)^2}{4}\right)$$

$$+ \omega^2\left(\frac{1}{\sqrt{n}}, \mathcal{F}_i\right) + \sum_{j=i+1}^{k} 6\omega^2\left(\frac{1}{\sqrt{n}}, \mathcal{F}_j\right) \cdot k^{1/2} \frac{A_i^2(n)}{A_j^2(n)}$$

$$\leq CA_i^2(n),$$

where $C$ is a constant not depending on $f$. Note that in the last inequality we use the fact that $\lambda^2 \exp(-\frac{(\lambda - 4(2k)^{1/2} - 2)^2}{4})$ is bounded as a function of $\lambda$. Hence

$$\sup_{f \in \mathcal{F}_i} E(\hat{T}^* - Tf)^2$$

$$= \max_{1\leq m\leq i}\left\{\sup_{f\in\mathcal{F}_m\setminus\mathcal{F}_{m-1}} E(\hat{T}^* - Tf)^2\right\}$$

$$\leq C\max_{1\leq m\leq i}\{A_m^2(n)\} = CA_i^2(n)$$

$$= C\left(\max_{1\leq m\leq i-1}\left\{\omega_+^2\left(\sqrt{\frac{\ln\gamma_{m,i,+}}{n}}, \mathcal{F}_m, \mathcal{F}_i\right)\right\} + \omega^2\left(\frac{1}{\sqrt{n}}, \mathcal{F}_i\right)\right). \quad \square$$

4.2. *Adaptation over nonnested parameter spaces.* Many common parameter spaces of interest such as Lipschitz spaces and Besov spaces are not nested. However, they often have nested structure in terms of the modulus of continuity. Theorem 4 can be generalized to such nonnested parameter spaces.

Let $\mathcal{F}_i$, $i = 1,\ldots,k$, be closed convex parameter spaces which are not necessarily nested. For any parameter set $\mathcal{F}$, let C.Hull($\mathcal{F}$) denote the convex hull of $\mathcal{F}$. We shall denote by $a(\varepsilon) \asymp b(\varepsilon)$ when $a(\varepsilon)/b(\varepsilon)$ is bounded away from 0 and $\infty$ as $\varepsilon \to 0+$. Suppose the parameter spaces $\mathcal{F}_i$ satisfy the following conditions on the modulus of continuity:

1. For $l \leq i$ and $m \leq j$, $\omega(\varepsilon, \mathcal{F}_l, \mathcal{F}_m) \leq C\omega(\varepsilon, \mathcal{F}_i, \mathcal{F}_j)$ for some constant $C > 0$.
2. For $2 \leq i \leq k$, $\omega(\varepsilon, \mathcal{G}_i) \asymp \omega(\varepsilon, \text{C.Hull}(\mathcal{G}_i))$ where $\mathcal{G}_i = \bigcup_{m=1}^i \mathcal{F}_m$.

Note that conditions 1 and 2 are trivially satisfied if $\mathcal{F}_i$ are nested.

As shown in [6], the minimax linear rate of convergence for estimating a linear functional $Tf$ over a parameter set $\mathcal{F}$ is determined by the modulus over its convex hull, $\omega(\frac{1}{\sqrt{n}}, \text{C.Hull}(\mathcal{F}))$. Conditions 1 and 2 together yield $\omega(\varepsilon, \mathcal{F}_i) \asymp \omega(\varepsilon, \text{C.Hull}(\mathcal{G}_i))$ and this consequently implies that for $1 \leq i \leq k$ there exists a rate optimal linear estimator $\hat{T}_i$ over $\mathcal{F}_i$ such that

(68) $$\sup_{f\in\bigcup_{m=1}^i \mathcal{F}_m} E(\hat{T}_i - Tf)^2 \leq C\omega^2\left(\frac{1}{\sqrt{n}}, \mathcal{F}_i\right).$$

Now define the quantities $\gamma_{i,j}$ and $\gamma_{i,j,+}$ as in (45), (46) and (47). Let $\hat{T}^*$ be defined the same as in Section 4.1 with the minimax rate optimal linear estimator $\hat{T}_i$ over $\mathcal{F}_i$ satisfying (68). Under conditions 1 and 2 above, the estimator $\hat{T}^*$ then achieves adaptation over the parameter spaces $\mathcal{F}_i$ with minimum cost. More precisely, we have the following.

THEOREM 5. *Let $\mathcal{F}_i$, $i = 1,\ldots,k$, be closed convex parameter spaces satisfying conditions 1 and 2 above and let the estimator $\hat{T}^*$ be given as above. Then $\hat{T}^*$ is optimally adaptive over $\mathcal{F}_i$ for $i = 1,\ldots,k$, that is,*

(69) $$\sup_{f\in\mathcal{F}_1} E(\hat{T}^* - Tf)^2 \leq C\omega^2\left(\frac{1}{\sqrt{n}}, \mathcal{F}_1\right),$$

*and for $2 \leq i \leq k$,*

$$\sup_{f \in \mathcal{F}_i} E(\hat{T}^* - Tf)^2$$
(70)
$$\leq C \left( \max_{1 \leq m \leq i-1} \left\{ \omega_+^2 \left( \sqrt{\frac{\ln \gamma_{m,i,+}}{n}}, \mathcal{F}_m, \mathcal{F}_i \right) \right\} + \omega^2 \left( \frac{1}{\sqrt{n}}, \mathcal{F}_i \right) \right).$$

The proof of Theorem 5 is essentially the same as that of Theorem 4. We omit the proof for reasons of space.

**5. Adaptation over infinitely many parameter spaces.** Section 4 gives a construction of adaptive estimators over collections of finitely many parameter spaces. In this section we shall further extend these results to a continuum of parameter spaces when the penalty of adaptation is always a logarithmic factor of the noise level. Well-known examples of such cases include estimating a function at a point over a collection of Lipschitz classes or Besov classes.

Let $\{\mathcal{F}_\lambda : \lambda \in \Lambda\}$ be a collection of closed convex parameter spaces and $T$ be a linear functional. Suppose that the following conditions hold for some constants $0 < c_1 \leq c_2 < \infty$ and $\varepsilon_0 > 0$:

C1. The index set $\Lambda$ is an ordered set with $\min(\Lambda) = \lambda_* \in \Lambda$, $\max(\Lambda) = \lambda^* \in \Lambda$, and $\mathcal{F}_{\lambda_1} \subset \mathcal{F}_{\lambda_2}$ if $\lambda_1 > \lambda_2$.
C2. For all $0 < \varepsilon \leq \varepsilon_0$ and all $\lambda \in \Lambda$, $c_1 \varepsilon^{r_\lambda} \leq \omega(\varepsilon, \mathcal{F}_\lambda) \leq c_2 \varepsilon^{r_\lambda}$ where $0 < r_{\lambda_2} < r_{\lambda_1} \leq 1$ if $\lambda_2 < \lambda_1$.
C3. For $\lambda_2 < \lambda_1$ and for all $0 < \varepsilon \leq \varepsilon_0$, $\omega(\varepsilon, \mathcal{F}_{\lambda_1}) < \omega(\varepsilon, \mathcal{F}_{\lambda_2})$ and $\omega_+(\varepsilon, \mathcal{F}_{\lambda_1}, \mathcal{F}_{\lambda_2}) \geq c_1 \omega(\varepsilon, \mathcal{F}_{\lambda_2})$.
C4. For any fixed $0 < \varepsilon \leq \varepsilon_0$ the set $\{\omega(\varepsilon, \mathcal{F}_\lambda) : \lambda \in \Lambda\}$ is compact.

Under these conditions it is clear from Theorem 1 that the minimum cost of adaptation is at least a logarithmic factor for any $\mathcal{F}_\lambda$ with $\lambda < \lambda^*$. We shall develop an adaptive procedure over the whole collection $\{\mathcal{F}_\lambda : \lambda \in \Lambda\}$ which attains the exact minimax rate of convergence over $\mathcal{F}_{\lambda^*}$ and attains the lower bounds given in Theorem 1 over any $\mathcal{F}_\lambda$ with $\lambda \in \Lambda$ and $\lambda < \lambda^*$.

The main idea behind the construction of the adaptive estimator is to first put down a finite grid of parameter spaces such that the modulus of continuity over each space on the grid is at least a fixed constant factor apart from the modulus for any other space on the grid; moreover, the modulus over any space in the collection $\{\mathcal{F}_\lambda : \lambda \in \Lambda\}$ is at most a fixed constant factor away from the modulus over one of the parameter spaces on this grid. We then use the techniques developed in Section 4.1 to construct a procedure which is adaptive over the finite grid. This procedure which is adaptive over the grid is then automatically adaptive over the whole collection. The construction of the grid is based on the following simple lemma.

LEMMA 5. *Let $\Omega$ be a compact subset of the positive half line $\mathbb{R}_+$ such that there exists an $\omega \in \Omega$ satisfying $\min(\Omega) < \omega \leq \frac{1}{2}\max(\Omega)$. Then there exists a unique finite sequence $\xi_1 < \xi_2 < \cdots < \xi_k$ with $\xi_i \in \Omega$, $\xi_1 = \min(\Omega)$ and $\xi_k = \max(\Omega)$ such that $\xi_{i+1} \geq 2\xi_i$ for all $2 \leq i \leq k-1$ and for any $\omega \in \Omega$ there exists $1 \leq i \leq k$ such that $\frac{1}{2}\xi_i < \omega \leq \xi_i$.*

The grid is constructed as follows. Set $\Omega_n = \{\omega(\sqrt{\frac{\log n}{n}}, \mathcal{F}_\lambda) : \lambda \in \Lambda\}$. Then for sufficiently large $n$ it follows from condition C4 that the set $\Omega_n$ is compact. If for all $\omega \in \Omega$ with $\omega > \min(\Omega)$, $\omega > \frac{1}{2}\max(\Omega)$, then set $k_n = 2$, $\xi_1 = \min(\Omega_n)$ and $\xi_{k_n} = \max(\Omega_n)$. Otherwise there is a sequence $\xi_1 < \xi_2 < \cdots < \xi_{k_n}$ in $\Omega_n$ satisfying the conditions given in Lemma 5. Let $\mathcal{F}_{\lambda_1} \subset \mathcal{F}_{\lambda_2} \subset \cdots \subset \mathcal{F}_{\lambda_{k_n}}$ be the corresponding closed convex parameter spaces with $\lambda_i \in \Lambda$ and $\xi_i = \omega(\sqrt{\frac{\log n}{n}}, \mathcal{F}_{\lambda_i})$. Note that it follows from the conditions $\lambda_1 = \lambda^* = \max(\Lambda)$, $\lambda_{k_n} = \lambda_* = \min(\Lambda)$ and $k_n \leq \log_2 n$ for large enough $n$. For convenience write $\mathcal{F}_i$ for $\mathcal{F}_{\lambda_i}$. This sequence of parameter spaces $\{\mathcal{F}_i : 1 \leq i \leq k_n\}$ forms a grid over the whole collection of parameter spaces $\{\mathcal{F}_\lambda : \lambda \in \Lambda\}$ such that for any $\lambda \in \Lambda$ with $\lambda < \lambda^*$, there exists $2 \leq i \leq k_n$ satisfying $\mathcal{F}_\lambda \subseteq \mathcal{F}_i$ and

$$\text{(71)} \qquad \frac{1}{2}\omega\left(\sqrt{\frac{\log n}{n}}, \mathcal{F}_i\right) \leq \omega\left(\sqrt{\frac{\log n}{n}}, \mathcal{F}_\lambda\right) \leq \omega\left(\sqrt{\frac{\log n}{n}}, \mathcal{F}_i\right).$$

We shall now turn to the construction of the adaptive estimator based on this grid. The construction is different but similar to the one given in Section 4.1. Let $\hat{T}_1$ be a linear estimator satisfying $\sup_{f \in \mathcal{F}_1} E(\hat{T}_1 - Tf)^2 \leq \omega^2(\frac{1}{\sqrt{n}}, \mathcal{F}_1)$. For $1 \leq i, j \leq k_n$ with $\max(i, j) \geq 2$ let $\hat{T}_{i,j}$ be the estimator satisfying (12)–(14) with $\mathcal{F} = \mathcal{F}_i$, $\mathcal{H} = \mathcal{F}_j$ and $V = \frac{1}{\log n}\omega^2(\sqrt{\frac{\log n}{n}}, \mathcal{F}_i, \mathcal{F}_j)$.

For $i < j$ the test between $\mathcal{F}_i$ and $\mathcal{F}_j$ is given by

$$\text{(72)} \quad I_{i,j} = 1\left(\hat{T}_{i,j} - \frac{11}{2}\omega\left(\sqrt{\frac{\log n}{n}}, \mathcal{F}_j\right) \leq \hat{T}_i \leq \hat{T}_{j,i} + \frac{11}{2}\omega\left(\sqrt{\frac{\log n}{n}}, \mathcal{F}_j\right)\right).$$

The test is in favor of $\mathcal{F}_i$ if $I_{i,j} = 1$. Our adaptive procedure is described sequentially in exactly the same way as given in Section 4.1. Formally the estimator $\hat{T}^*$ can be written as

$$\text{(73)} \qquad \hat{T}^* = \sum_{i=1}^{k_n} \left(1 - \prod_{m<i}\prod_{j>m} I_{m,j}\right)\left(\prod_{j>i} I_{i,j}\right)\hat{T}_i.$$

The following theorem shows that this estimator is optimally adaptive over the whole collection of parameter spaces $\{\mathcal{F}_\lambda : \lambda \in \Lambda\}$.

THEOREM 6. *Let $\{\mathcal{F}_\lambda : \lambda \in \Lambda\}$ be a collection of nested closed convex parameter spaces and $T$ be a linear functional. Suppose that conditions C1–C4 hold. Then the estimator $\hat{T}^*$ defined in (73) is optimally adaptive over $\mathcal{F}_\lambda$ for all $\lambda \in \Lambda$. More specifically, there exists a constant $C > 0$ such that*

$$(74) \qquad \sup_{f \in \mathcal{F}_{\lambda^*}} E(\hat{T}^* - Tf)^2 \leq C\omega^2\left(\frac{1}{\sqrt{n}}, \mathcal{F}_{\lambda^*}\right),$$

*and for all $\lambda \in \Lambda$ and $\lambda < \lambda^*$*

$$(75) \qquad \sup_{f \in \mathcal{F}_\lambda} E(\hat{T}^* - Tf)^2 \leq C\omega^2\left(\sqrt{\frac{\log n}{n}}, \mathcal{F}_\lambda\right).$$

REMARK. The structural conditions C1–C4 are used to keep track of a growing number of between-class moduli. These conditions seem to be necessary for developing an adaptation theory over infinitely many parameter spaces. The completely general setting is difficult because it is possible that the penalty for adaptation varies from space to space in a very complicated way from no penalty to a logarithmic factor to an algebraic factor.

The proof of Theorem 6 is similar to that of Theorem 4. It relies on the analysis of the tests $I_{i,j}$. For $f \in \mathcal{F}_i \setminus \mathcal{F}_{i-1}$, the main analysis is concerned with the cases where $\hat{T}_j$ with $j < i$ is used and where $\hat{T}_j$ with $j > i$ is used. Lemma 6 shows that the chance of using $\hat{T}_j$ with $j > i$ is small. Lemma 7 shows that the chance of using $\hat{T}_j$ with $j < i$ is small whenever the bias of $\hat{T}_j$ is large. Before presenting these technical results we first collect some useful bounds on the expectations and variances of $\hat{T}_i - \hat{T}_{i,j}$ and $\hat{T}_i - \hat{T}_{j,i}$.

For $j \geq 2$ set

$$(76) \qquad b_j = \frac{3}{2}\omega\left(\sqrt{\frac{\log n}{n}}, \mathcal{F}_j\right) \quad \text{and} \quad v_j = \frac{4}{\log n}\omega^2\left(\sqrt{\frac{\log n}{n}}, \mathcal{F}_j\right)$$

and write $\hat{T}_j$ for $\hat{T}_{j,j}$. Note that if $f \in \mathcal{F}_i$, then for $j \geq i$,

$$(77) \quad E(\hat{T}_i - \hat{T}_{i,j}) \geq -\omega\left(\sqrt{\frac{\log n}{n}}, \mathcal{F}_i\right) - \frac{1}{2}\omega\left(\sqrt{\frac{\log n}{n}}, \mathcal{F}_i, \mathcal{F}_j\right) \geq -b_j,$$

$$(78) \quad E(\hat{T}_i - \hat{T}_{j,i}) \leq \omega\left(\sqrt{\frac{\log n}{n}}, \mathcal{F}_i\right) + \frac{1}{2}\omega\left(\sqrt{\frac{\log n}{n}}, \mathcal{F}_j, \mathcal{F}_i\right) \leq b_j.$$

For the variances, note that if $f \in \mathcal{F}_1$ and $j \geq 2$,

$$\operatorname{Var}(\hat{T}_1 - \hat{T}_{1,j}) \leq 2\left(\omega^2\left(\frac{1}{\sqrt{n}}, \mathcal{F}_1\right) + \frac{1}{\log n}\omega^2\left(\sqrt{\frac{\log n}{n}}, \mathcal{F}_1, \mathcal{F}_j\right)\right)$$

(79)
$$\leq 2\left(\omega^2\left(\frac{1}{\sqrt{n}}, \mathcal{F}_1\right) + \frac{1}{\log n}\omega^2\left(\sqrt{\frac{\log n}{n}}, \mathcal{F}_j\right)\right),$$

(80)
$$\operatorname{Var}(\hat{T}_1 - \hat{T}_{j,1}) \leq 2\left(\omega^2\left(\frac{1}{\sqrt{n}}, \mathcal{F}_1\right) + \frac{1}{\log n}\omega^2\left(\sqrt{\frac{\log n}{n}}, \mathcal{F}_j, \mathcal{F}_1\right)\right)$$
$$\leq 2\left(\omega^2\left(\frac{1}{\sqrt{n}}, \mathcal{F}_1\right) + \frac{1}{\log n}\omega^2\left(\sqrt{\frac{\log n}{n}}, \mathcal{F}_j\right)\right),$$

and if $f \in \mathcal{F}_i$ with $2 \leq i \leq j$,

(81)
$$\operatorname{Var}(\hat{T}_i - \hat{T}_{i,j}) \leq 2\left(\frac{1}{\log n}\omega^2\left(\sqrt{\frac{\log n}{n}}, \mathcal{F}_i\right) + \frac{1}{\log n}\omega^2\left(\sqrt{\frac{\log n}{n}}, \mathcal{F}_i, \mathcal{F}_j\right)\right)$$
$$\leq v_j,$$

(82)
$$\operatorname{Var}(\hat{T}_i - \hat{T}_{j,i}) \leq 2\left(\frac{1}{\log n}\omega^2\left(\sqrt{\frac{\log n}{n}}, \mathcal{F}_i\right) + \frac{1}{\log n}\omega^2\left(\sqrt{\frac{\log n}{n}}, \mathcal{F}_j, \mathcal{F}_i\right)\right)$$
$$\leq v_j.$$

LEMMA 6. *If* $f \in \mathcal{F}_1$, *then*

(83) $$P(\hat{T}^* = \hat{T}_2) \leq 4\exp\left(-2\frac{\omega^2(\sqrt{\log n/n}, \mathcal{F}_2)}{\omega^2(1/\sqrt{n}, \mathcal{F}_1)}\right) + 2k_n n^{-2},$$

*and for* $j \geq 3$,

(84) $$P(\hat{T}^* = \hat{T}_j) \leq \min_{i \leq m \leq j-1} P(I_m = 0) \leq 2k_n n^{-2}.$$

*If* $f \in \mathcal{F}_i$ *with* $i \geq 2$, *then for* $j > i$,

(85) $$P(\hat{T}^* = \hat{T}_j) \leq 2k_n n^{-2}.$$

*In particular,* (85) *holds for* $f \in \mathcal{F}_i \setminus \mathcal{F}_{i-1}$.

PROOF. First note that (84) follows from (85) so we need only prove (83) and (85). We first prove (85). It follows from (77), (78), (81) and (82) that for $f \in \mathcal{F}_i$ and $2 \leq i < j \leq k_n$,

$$P(I_{i,j} = 0) \leq P\left(\hat{T}_i - \hat{T}_{i,j} \leq -\frac{11}{2}\omega\left(\sqrt{\frac{\log n}{n}}, \mathcal{F}_j\right)\right)$$

$$+ P\left(\hat{T}_i - \hat{T}_{j,i} \geq \frac{11}{2}\omega\left(\sqrt{\frac{\log n}{n}}, \mathcal{F}_j\right)\right)$$

$$\leq \exp\left(-\frac{((11/2)\omega(\sqrt{\log n/n}, \mathcal{F}_j) - b_j)^2}{2v_j}\right)$$

$$+ \exp\left(-\frac{((11/2)\omega(\sqrt{\log n/n}, \mathcal{F}_j) - b_j)^2}{2v_j}\right)$$

$$\leq 2\exp\left(-\frac{(11/2 - 3/2)^2\omega^2(\sqrt{\log n/n}, \mathcal{F}_j)}{(8/\log n)\omega^2(\sqrt{\log n/n}, \mathcal{F}_j)}\right) = 2n^{-2}.$$

Set $I_i = \prod_{j>i} I_{i,j}$. Then $P(I_m = 0) \leq \sum_{l=m+1}^{k_n} P(I_{m,l} = 0) \leq 2k_n n^{-2}$ and hence

(86) $$P(\hat{T}^* = \hat{T}_j) \leq \min_{i \leq m \leq j-1} P(I_m = 0) \leq 2k_n n^{-2}.$$

This proves (85). Now assume that $f \in \mathcal{F}_1$. Then for $j \geq 2$,

$$P(I_{1,j} = 0) \leq P\left(\hat{T}_1 - \hat{T}_{1,j} \leq -\frac{11}{2}\omega\left(\sqrt{\frac{\log n}{n}}, \mathcal{F}_j\right)\right)$$

$$+ P\left(\hat{T}_1 - \hat{T}_{j,1} \geq \frac{11}{2}\omega\left(\sqrt{\frac{\log n}{n}}, \mathcal{F}_j\right)\right)$$

$$\leq 2\exp\left(-\frac{(11/2 - 3/2)^2\omega^2(\sqrt{\log n/n}, \mathcal{F}_j)}{4\omega^2(1/\sqrt{n}, \mathcal{F}_1) + (4/\log n)\omega^2(\sqrt{\log n/n}, \mathcal{F}_j)}\right).$$

We consider two cases. First if $\omega^2(\frac{1}{\sqrt{n}}, \mathcal{F}_1) \leq \frac{1}{\log n}\omega^2(\sqrt{\frac{\log n}{n}}, \mathcal{F}_2)$, then

$$P(I_{1,j} = 0) \leq 2\exp\left(-\frac{(11/2 - 3/2)^2\omega^2(\sqrt{\log n/n}, \mathcal{F}_j)}{(8/\log n)\omega^2(\sqrt{\log n/n}, \mathcal{F}_j)}\right) = 2n^{-2}$$

and hence $P(\hat{T}^* = \hat{T}_2) \leq 2k_n n^{-2}$ and (83) follows.

Now suppose that $\omega^2(\frac{1}{\sqrt{n}}, \mathcal{F}_1) \geq \frac{1}{\log n}\omega^2(\sqrt{\frac{\log n}{n}}, \mathcal{F}_2)$. Let $j = j_*$ be the largest integer such that $\omega^2(\frac{1}{\sqrt{n}}, \mathcal{F}_1) \geq \frac{1}{\log n}\omega^2(\sqrt{\frac{\log n}{n}}, \mathcal{F}_j)$. Then for $j_* + 1 \leq j \leq k_n$ it is easy to see that $P(I_{1,j} = 0) \leq 2n^{-2}$. For $2 \leq j \leq j_*$,

$$P(I_{1,j} = 0) \leq 2\exp\left(-\frac{(11/2 - 3/2)^2\omega^2(\sqrt{\log n/n}, \mathcal{F}_j)}{4\omega^2(1/\sqrt{n}, \mathcal{F}_1) + (4/\log n)\omega^2(\sqrt{\log n/n}, \mathcal{F}_j)}\right)$$

$$\leq 2\exp\left(-2\frac{\omega^2(\sqrt{\log n/n}, \mathcal{F}_j)}{\omega^2(1/\sqrt{n}, \mathcal{F}_1)}\right).$$

Note that by the construction of the grid of parameter spaces $\{\mathcal{F}_i : 1 \leq i \leq k_n\}$ for $j \geq 2$,

$$\omega^2\left(\sqrt{\frac{\log n}{n}}, \mathcal{F}_j\right) \geq 2^{2j-4}\omega^2\left(\sqrt{\frac{\log n}{n}}, \mathcal{F}_2\right).$$

Hence,

$$P(\hat{T}^* = \hat{T}_2) \leq P(I_1 = 0) \leq \sum_{j=2}^{k_n} P(I_{1,j} = 0)$$

$$\leq 2\sum_{j=2}^{j_*} \exp\left(-2\frac{\omega^2(\sqrt{\log n/n}, \mathcal{F}_j)}{\omega^2(1/\sqrt{n}, \mathcal{F}_1)}\right) + 2\sum_{j=j_*+1}^{k_n} n^{-2}$$

$$\leq 2\sum_{j=2}^{j_*} \exp\left(-2^{2j-4} \cdot 2\frac{\omega^2(\sqrt{\log n/n}, \mathcal{F}_2)}{\omega^2(1/\sqrt{n}, \mathcal{F}_1)}\right) + 2k_n n^{-2}$$

$$\leq 4\exp\left(-2\frac{\omega^2(\sqrt{\log n/n}, \mathcal{F}_2)}{\omega^2(1/\sqrt{n}, \mathcal{F}_1)}\right) + 2k_n n^{-2}$$

and once again (83) follows. $\square$

LEMMA 7. *Suppose $f \in \mathcal{F}_i \setminus \mathcal{F}_{i-1}$ and $j < i$. If $|E\hat{T}_j - Tf| \geq \beta\omega(\sqrt{\frac{\log n}{n}}, \mathcal{F}_i)$ for some $\beta > 6$, then*

(87) $$P(\hat{T}^* = \hat{T}_j) \leq n^{-(\beta-6)^2/8}.$$

PROOF. We shall only consider the case when $E\hat{T}_j - Tf \geq \beta\omega(\sqrt{\frac{\log n}{n}}, \mathcal{F}_i)$ since the case of $E\hat{T}_j - Tf \leq -\beta\omega(\sqrt{\frac{\log n}{n}}, \mathcal{F}_i)$ can be handled similarly. Let $f \in \mathcal{F}_i \setminus \mathcal{F}_{i-1}$. Then $P(\hat{T}^* = T_j) \leq P(I_{j,i} = 1) \leq P(\hat{T}_j - \hat{T}_{i,j} - \frac{11}{2}\omega(\sqrt{\frac{\log n}{n}}, \mathcal{F}_i) \leq 0)$. Note that

$$E\left(\hat{T}_j - \hat{T}_{i,j} - \frac{11}{2}\omega\left(\sqrt{\frac{\log n}{n}}, \mathcal{F}_i\right)\right)$$

$$= E(\hat{T}_j - Tf) - E(\hat{T}_{i,j} - Tf) - \frac{11}{2}\omega\left(\sqrt{\frac{\log n}{n}}, \mathcal{F}_i\right)$$

$$\geq \beta\omega\left(\sqrt{\frac{\log n}{n}}, \mathcal{F}_i\right) - \frac{1}{2}\omega\left(\sqrt{\frac{\log n}{n}}, \mathcal{F}_i, \mathcal{F}_j\right) - \frac{11}{2}\omega\left(\sqrt{\frac{\log n}{n}}, \mathcal{F}_i\right)$$

$$\geq (\beta - 6)\omega\left(\sqrt{\frac{\log n}{n}}, \mathcal{F}_i\right).$$

Note that $\operatorname{Var}(\hat{T}_j - \hat{T}_{i,j}) \leq \frac{4}{\log n}\omega^2(\sqrt{\frac{\log n}{n}}, \mathcal{F}_i)$ and hence

$$P(\hat{T}^* = T_j) \leq \exp\left(-\frac{(\beta-6)^2}{2} \frac{\omega^2(\sqrt{\log n/n}, \mathcal{F}_i)}{(4/\log n)\omega^2(\sqrt{\log n/n}, \mathcal{F}_i)}\right)$$
$$\leq n^{-(\beta-6)^2/8}. \qquad \square$$

We are now ready to prove Theorem 6 using the above technical results. We first show that the estimator $\hat{T}^*$ given in (73) has the desired adaptation properties over the grid $\{\mathcal{F}_i : 1 \leq i \leq k_n\}$ and then show that $\hat{T}^*$ is in fact adaptive over the whole collection of parameter spaces $\{\mathcal{F}_\lambda : \lambda \in \Lambda\}$.

PROOF OF THEOREM 6. The proof is broken into three steps. In each step it is important to note that $k_n \leq \log_2 n$ and $\omega(\sqrt{\frac{\log n}{n}}, \mathcal{F}_j) \leq \frac{1}{2}\omega(\sqrt{\frac{\log n}{n}}, \mathcal{F}_{j+1})$.

*Step* 1. We begin by showing that the estimator $\hat{T}^*$ attains the exact minimax rate over $\mathcal{F}_1 = \mathcal{F}_{\lambda^*}$. We shall only consider the case $\omega^2(\frac{1}{\sqrt{n}}, \mathcal{F}_1) \geq \frac{1}{\log n}\omega^2(\sqrt{\frac{\log n}{n}}, \mathcal{F}_2)$. When $\omega^2(\frac{1}{\sqrt{n}}, \mathcal{F}_1) < \frac{1}{\log n}\omega^2(\sqrt{\frac{\log n}{n}}, \mathcal{F}_2)$ the proof is easier. Note that

$$\sup_{f \in \mathcal{F}_1} E(\hat{T}^* - Tf)^2 \leq \sup_{f \in \mathcal{F}_1} E(\hat{T}_1 - Tf)^2$$
$$+ \sup_{f \in \mathcal{F}_1} (E|\hat{T}_2 - Tf|^4)^{1/2} \cdot (P(\hat{T}^* = \hat{T}_2))^{1/2}$$
$$+ \sum_{j=3}^{k_n} \sup_{f \in \mathcal{F}_1} (E|\hat{T}_j - Tf|^4)^{1/2} \cdot (P(\hat{T}^* = \hat{T}_j))^{1/2}$$
$$\leq \omega^2\left(\frac{1}{\sqrt{n}}, \mathcal{F}_1\right)$$
$$+ C\omega^2\left(\sqrt{\frac{\log n}{n}}, \mathcal{F}_2\right) \cdot 2\exp\left(-\frac{\omega^2(\sqrt{\log n/n}, \mathcal{F}_2)}{\omega^2(1/\sqrt{n}, \mathcal{F}_1)}\right)$$
$$+ C\sum_{j=2}^{k_n} \omega^2\left(\sqrt{\frac{\log n}{n}}, \mathcal{F}_j\right) \cdot (2k_n)^{1/2} n^{-1}.$$

Now note that
$$\omega^2\left(\sqrt{\frac{\log n}{n}}, \mathcal{F}_2\right) \exp\left(-\frac{\omega^2(\sqrt{\log n/n}, \mathcal{F}_2)}{\omega^2(1/\sqrt{n}, \mathcal{F}_1)}\right)$$
$$\leq \omega^2\left(\sqrt{\frac{\log n}{n}}, \mathcal{F}_1\right) \cdot \frac{\omega^2(\sqrt{\log n/n}, \mathcal{F}_2)}{\omega^2(1/\sqrt{n}, \mathcal{F}_1)} \exp\left(-\frac{\omega^2(\sqrt{\log n/n}, \mathcal{F}_2)}{\omega^2(1/\sqrt{n}, \mathcal{F}_1)}\right)$$

$$\leq \omega^2\left(\sqrt{\frac{\log n}{n}}, \mathcal{F}_1\right) \sup_{x \geq 2} x e^{-x} \leq \frac{1}{2}\omega^2\left(\sqrt{\frac{\log n}{n}}, \mathcal{F}_1\right)$$

and

$$\sum_{j=2}^{k_n} \omega^2\left(\sqrt{\frac{\log n}{n}}, \mathcal{F}_j\right) \cdot (2k_n)^{1/2} n^{-1} \leq 2\omega^2\left(\sqrt{\frac{\log n}{n}}, \mathcal{F}_{k_n}\right) \cdot (2k_n)^{1/2} n^{-1}$$
$$= o(n^{-1}).$$

Hence $\sup_{f \in \mathcal{F}_1} E(\hat{T}^* - Tf)^2 \leq C\omega^2(\frac{1}{\sqrt{n}}, \mathcal{F}_1)$ for some absolute constant $C > 0$.

*Step* 2. Now consider $i \geq 2$. Let $f \in \mathcal{F}_i \setminus \mathcal{F}_{i-1}$ for some $i \geq 2$. Set

$$\mathcal{J}_1 = \left\{ j < i : |E\hat{T}_j - Tf| \leq 7\omega\left(\sqrt{\frac{\log n}{n}}, \mathcal{F}_i\right) \right\},$$

$$\mathcal{J}_2 = \left\{ j < i : |E\hat{T}_j - Tf| > 7\omega\left(\sqrt{\frac{\log n}{n}}, \mathcal{F}_i\right) \right\}.$$

Then for $f \in \mathcal{F}_i \setminus \mathcal{F}_{i-1}$ with $i > 1$ we have

$$E(\hat{T}^* - Tf)^2 = \sum_{j=1}^{k_n} E\{(\hat{T}_j - Tf)^2 1(\hat{T}^* = \hat{T}_j)\}$$
$$\leq \sum_{j \in \mathcal{J}_1} E\{(\hat{T}_j - Tf)^2 1(\hat{T}^* = \hat{T}_j)\}$$
$$+ \sum_{j \in \mathcal{J}_2} (E|\hat{T}_j - Tf|^4)^{1/2} (P(\hat{T}^* = \hat{T}_j))^{1/2}$$
$$+ E(\hat{T}_i - Tf)^2 + \sum_{j=i+1}^{k_n} (E|\hat{T}_j - Tf|^4)^{1/2} (P(\hat{T}^* = \hat{T}_j))^{1/2}$$
$$\equiv S_1 + S_2 + S_3 + S_4.$$

We bound the four terms separately. First consider $S_1$. Note that

$$S_1 = \sum_{j \in \mathcal{J}_1} E\{(\hat{T}_j - Tf)^2 1(\hat{T}^* = \hat{T}_j)\}$$
$$\leq 2 \sum_{j \in \mathcal{J}_1} E\{[(\hat{T}_j - E\hat{T}_j)^2 + (E\hat{T}_j - Tf)^2] 1(\hat{T}^* = \hat{T}_j)\}$$
$$\leq 2 \sum_{j \in \mathcal{J}_1} \text{Var}(\hat{T}_j) + 2 \sum_{j \in \mathcal{J}_1} (E\hat{T}_j - Tf)^2 P(\hat{T}^* = \hat{T}_j)$$

$$\leq 2 \sum_{j=1}^{i-1} \frac{1}{\log n} \omega^2 \left( \sqrt{\frac{\log n}{n}}, \mathcal{F}_j \right) + 98 \omega^2 \left( \sqrt{\frac{\log n}{n}}, \mathcal{F}_i \right)$$

$$\leq 100 \omega^2 \left( \sqrt{\frac{\log n}{n}}, \mathcal{F}_i \right).$$

Now consider $S_2$. If $j \in \mathcal{J}_2$, then $|E\hat{T}_j - Tf| = \beta_j \omega(\sqrt{\frac{\log n}{n}}, \mathcal{F}_i)$ for some $\beta_j > 7$. Hence by (30),

$$(E|\hat{T}_j - Tf|^4)^{1/2} \leq 3 \operatorname{Var}(\hat{T}_j) + 3\beta_j^2 \omega^2 \left( \sqrt{\frac{\log n}{n}}, \mathcal{F}_i \right)$$

$$\leq 4\beta_j^2 \omega^2 \left( \sqrt{\frac{\log n}{n}}, \mathcal{F}_i \right).$$

It follows from Lemma 7 that $(P(\hat{T}^* = T_j))^{1/2} \leq n^{-(\beta_j - 6)^2/16}$. Hence

$$S_2 = \sum_{j \in \mathcal{J}_2} (E|\hat{T}_j - Tf|^4)^{1/2} (P(\hat{T}^* = \hat{T}_j))^{1/2}$$

$$\leq \omega^2 \left( \sqrt{\frac{\log n}{n}}, \mathcal{F}_i \right) \cdot \sum_{j \in \mathcal{J}_2} 4\beta_j^2 n^{-(\beta_j - 6)^2/16}$$

$$\leq \omega^2 \left( \sqrt{\frac{\log n}{n}}, \mathcal{F}_i \right) \cdot 4 k_n n^{-1/16} \sup_{x \geq 7} x^2 n^{-((x-6)^2 - 1)/16}$$

$$= o\left( \omega^2 \left( \sqrt{\frac{\log n}{n}}, \mathcal{F}_i \right) \right).$$

For $S_3$ it is clear from the construction of $\hat{T}_i = \hat{T}_{i,i}$ that

$$S_3 = E(\hat{T}_i - Tf)^2 \leq 2\omega^2 \left( \sqrt{\frac{\log n}{n}}, \mathcal{F}_i \right).$$

Finally for $S_4$ it follows from (30) and Lemma 6 that

$$S_4 = \sum_{j=i+1}^{k_n} (E|\hat{T}_j - Tf|^4)^{1/2} (P(\hat{T}^* = \hat{T}_j))^{1/2}$$

$$\leq \sum_{j=i+1}^{k_n} 3 \left( \frac{1}{4} \omega^2 \left( \sqrt{\frac{\log n}{n}}, \mathcal{F}_j \right) + \frac{1}{\log n} \omega^2 \left( \sqrt{\frac{\log n}{n}}, \mathcal{F}_j \right) \right) \cdot (2k_n)^{1/2} n^{-1}$$

$$= o(n^{-1}).$$

ON ADAPTIVE ESTIMATION OF LINEAR FUNCTIONALS   33

Putting the four terms together we have that, for $f \in \mathcal{F}_i \setminus \mathcal{F}_{i-1}$ with $i \geq 2$,

$$E(\hat{T}^* - Tf)^2 \leq S_1 + S_2 + S_3 + S_4 \leq C\omega^2\left(\sqrt{\frac{\log n}{n}}, \mathcal{F}_i\right),$$

where $C$ is an absolute constant not depending on $f$, $n$, $k_n$ and $i$. Hence for all $2 \leq i \leq k_n$,

$$\sup_{f \in \mathcal{F}_i} E(\hat{T}^* - Tf)^2 = \max_{1 \leq m \leq i}\left\{\sup_{f \in \mathcal{F}_m \setminus \mathcal{F}_{m-1}} E(\hat{T}^* - Tf)^2\right\}$$

$$\leq C \max_{1 \leq m \leq i}\left\{\omega^2\left(\sqrt{\frac{\log n}{n}}, \mathcal{F}_m\right)\right\} = C\omega^2\left(\sqrt{\frac{\log n}{n}}, \mathcal{F}_i\right).$$

*Step* 3. Steps 1 and 2 show that the estimator $\hat{T}^*$ is adaptive over the grid $\{\mathcal{F}_i : 1 \leq i \leq k_n\}$. It is now easy to show that $\hat{T}^*$ is in fact adaptive over the collection $\{\mathcal{F}_\lambda : \lambda \in \Lambda\}$. Note that for any $\lambda \in \Lambda$ with $\lambda < \lambda^*$, by the construction of the grid $\{\mathcal{F}_i : 1 \leq i \leq k_n\}$, there exists $2 \leq i \leq k_n$ such that $\mathcal{F}_\lambda \subseteq \mathcal{F}_i$ with

$$\frac{1}{4}\omega^2\left(\sqrt{\frac{\log n}{n}}, \mathcal{F}_i\right) \leq \omega^2\left(\sqrt{\frac{\log n}{n}}, \mathcal{F}_\lambda\right) \leq \omega^2\left(\sqrt{\frac{\log n}{n}}, \mathcal{F}_i\right).$$

Hence

$$\sup_{f \in \mathcal{F}_\lambda} E(\hat{T}^* - Tf)^2 \leq \sup_{f \in \mathcal{F}_i} E(\hat{T}^* - Tf)^2$$

$$\leq C\omega^2\left(\sqrt{\frac{\log n}{n}}, \mathcal{F}_i\right) \leq 4C\omega^2\left(\sqrt{\frac{\log n}{n}}, \mathcal{F}_\lambda\right)$$

and the theorem is proved. $\square$

REMARK. Similarly to the case of finitely many parameter spaces, the results given above for infinitely many nested spaces can be extended in a straightforward way to nonnested parameter spaces when the moduli of continuity have nested structure under conditions similar to those given in Section 4.2.

Department of Statistics
The Wharton School
University of Pennsylvania
Philadelphia, Pennsylvania 19104-6340
USA
e-mail: tcai@wharton.upenn.edu